\pdfoutput=1
\documentclass[11pt,fullpage]{article} 
\usepackage[centertags]{amsmath}
\usepackage{amsfonts}
\usepackage{amssymb}
\usepackage{amsthm}
\usepackage{eepic}
\usepackage{newlfont}
\usepackage{graphicx}

\theoremstyle{plain}
\newtheorem{thm}{Theorem}[section]
\newcommand{\BTHM}{\begin{thm}} \newcommand{\ETHM}{\end{thm}}
\newtheorem{cor}[thm]{Corollary}
\newcommand{\BCR}{\begin{cor}} \newcommand{\ECR}{\end{cor}}
\newtheorem{lem}[thm]{Lemma}
\newcommand{\BL}{\begin{lem}}   \newcommand{\EL}{\end{lem}}
\newtheorem{clm}[thm]{Claim}
\newcommand{\BCM}{\begin{clm}}   \newcommand{\ECM}{\end{clm}}
\newtheorem{prop}[thm]{Proposition}
\newcommand{\BP}{\begin{prop}}   \newcommand{\EP}{\end{prop}}
\newtheorem{assm}[thm]{Assumption}
\newcommand{\BASM}{\begin{assm}}   \newcommand{\EASM}{\end{assm}}

\theoremstyle{definition}
\newtheorem{defn}{Definition}[section]
\newcommand{\BD}{\begin{defn}}   \newcommand{\ED}{\end{defn}}
\newtheorem{con}[thm]{Conjecture}
\newcommand{\BCONJ}{\begin{con}}   \newcommand{\ECONJ}{\end{con}}

\theoremstyle{definition}
\newtheorem{problem}[thm]{Problem}
\newcommand{\BPR}{\begin{problem}}   \newcommand{\EPR}{\end{problem}}

\newenvironment{rem}{\noindent{\bf Remark:~~}}{}
\newcommand{\BREM}{\begin{rem}} \newcommand{\EREM}{\end{rem}}
\newenvironment{discussion}{\noindent{\bf Discussion:~~\\}}{}
\newcommand{\BDIS}{\begin{discussion}} \newcommand{\EDIS}{\end{discussion}}

\numberwithin{equation}{section}

\def\blackslug
     {\hbox{\hskip 1pt\vrule width 8pt height 8pt depth 1.5pt\hskip 1pt}}
\def\qed{\quad\blackslug\lower 8.5pt\null\par}

\newtheorem{exmp}[thm]{Example}
\newcommand{\BEX}{\begin{exmp}} \newcommand{\EEX}{\end{exmp}}
\newcommand{\BF}{\begin{fact}}   \newcommand{\EF}{\end{fact}}
\newcommand{\Bcr}{\begin{techcorr}}
\newcommand{\Ecr}{\end{techcorr}}
\newcommand{\BDS}{\begin{description}}
\newcommand{\EDS}{\end{description}}
\newcommand{\BE}{\begin{enumerate}}
\newcommand{\EE}{\end{enumerate}}
\newcommand{\BI}{\begin{itemize}}
\newcommand{\EI}{\end{itemize}}
\renewenvironment{proof}{\noindent{\bf Proof:~~}}{\qed}
\newcommand{\BPF}{\begin{proof}}
\newcommand{\EPF}{\end{proof}}

\newcommand{\BB}{\begin{enumerate}}
\newcommand{\EB}{\end{enumerate}}

\begin{document}

\title{On the Profile of Multiplicities of Complete Subgraphs}

\author{Uriel Feige \qquad \qquad Anne Kenyon \qquad \qquad Shimon Kogan \\ \\
  Department of Computer Science and Applied Mathematics \\
          Weizmann Institue, Rehovot 76100, Israel \\ \\
           uriel.feige@weizmann.ac.il \\
           annemkenyon@gmail.com \\
           shimon.kogan@weizmann.ac.il
          }



\maketitle

\begin{abstract}
Let $G$ be a $2$-coloring of a complete graph on $n$ vertices, for sufficiently large $n$. We prove that $G$ contains at least $n^{(\frac{1}{4} - o(1))\log n}$ monochromatic complete subgraphs of size $r$, where
\[
0.3\log n <  r   < 0.7\log n.
\]
The previously known lower bound on the total number of monochromatic complete subgraphs, due to Sz{\'{e}}kely~\cite{DBLP:journals/combinatorica/Szekely84a}, was $n^{0.1576\log n}$.
We also prove that $G$ contains at least $n^{\frac{1}{7} \log n} $ monochromatic complete subgraphs of size $\frac{1}{2}\log n$.

If furthermore one assumes that the largest monochromatic complete subgraph in $G$ is of size $(\frac{1}{2} + o(1))\log n$ (it is a well known open question whether such graphs exist), then for every constant $0 \le c \le \frac{1}{2}$ we determine (up to low order terms) the number of monochromatic complete subgraphs of size $c \log n$. We do so by proving a lower bound that matches (up to low order terms) a previous upper bound of Sz{\'{e}}kely~\cite{DBLP:journals/combinatorica/Szekely84a}. For example, the number of monochromatic complete subgraphs of size $\frac{1}{2} \log n$ is
$n^{\frac{1}{8}(4 - \log e \pm o(1))\log n} \simeq n^{0.32 \log n}$.

\end{abstract}




\section{Introduction}
\label{sec:introduction}

The ``classic" diagonal Ramsey number question asks the following: what is the
minimum $n$ such that all $2$-colored (edge colored) complete graphs of size $n$
contain a monochromatic complete subgraph of size $t$?
Letting $\alpha$ be the size of the largest
monochromatic complete subgraph, the diagonal Ramsey number question can be rephrased
as: what is the minimum $\alpha$ over all possible $2$-colorings of the complete graph
on $n$ vertices?

These questions are the basis for the field of Ramsey Theory.
Known bounds say that all graphs of size
$n$ contain a monochromatic complete subgraph of size at least $\frac{\log n}{2}$ and that there are
graphs with a maximum monochromatic complete subgraph of size at most $2 \log n$. These results
date back to 1935 (Erd{\"o}s and Szekeres \cite{erdios1935combinatorial})
and 1947 (Erd{\"o}s \cite{erdos1947some}), respectively.
There have since been improvements to these bounds, but only to the lower order
terms.

In this paper we turn our
attention to a related question that is referred to as ``Ramsey Multiplicity":
what is the minimum number of monochromatic complete subgraphs of size $t$ in a $2$-colored complete
graph? The classic Ramsey problem is a special case of the Ramsey multiplicity question, in the sense that it asks to determine the largest value $t$ (as a function of $n$) for which the Ramsey multiplicity is guaranteed to be nonzero. Hence beyond its intrinsic interest, progress on the Ramsey Multiplicity question may potentially lead to progress on the classic Ramsey problem. Our interest in this work is mainly in the case that $t = \Theta(\log n)$, which is the relevant range for the classic Ramsey problem. There has also been previous work for the case of constant $t$ (see the Section~\ref{sec:relatedwork} for more details).

\subsection{Main results}

The question of the total number of monochromatic complete subgraphs is addressed by Szekely~\cite{DBLP:journals/combinatorica/Szekely84a}, who shows that for large enough $t$, in any $2$-coloring of a complete graph on $2^t$ vertices there are at
least $2^{0.1576t^2}$ monochromatic complete subgraphs. We improve this result by proving the following in Section \ref{central_results1}.

\BTHM
\label{thm:improveSzekely}
Let $G$ be a $2$-coloring of a complete graph on $2^t$ vertices. Then for any large enough $t$, graph $G$ contains at least
$2^{\left(\frac{1}{4}-o(1)\right)t^2}$ monochromatic complete subgraphs of size $r$ where
\[
0.3t < r < 0.7t
\]
Equivalently, every 2-coloring of a complete graph on $n$ vertices gives at least roughly $n^{\frac{1}{4}\log n}$ monochromatic subgraphs of size in the range $[0.3 \log n, 0.7 \log n]$.
\ETHM

Random graphs provide the known upper bound on the number of monochromatic complete subgraphs (see~\cite{DBLP:journals/combinatorica/Szekely84a} for details).

\BTHM
For all large enough $t$, there is a $2$-coloring of a complete graph on $2^t$ vertices with at most
$2^{(\frac{1}{2}+o(1))t^2}$ monochromatic complete subgraphs.
\ETHM

As noted above, it is a long standing open question whether for every $\epsilon > 0$ there are 2-colorings of the complete graph of size $n$ that do not induce a monochromatic complete subgraph of size $(\frac{1}{2}+\epsilon)\log n$. Not wishing to carry $\epsilon$ in our notation, we refer to the upper bound on the size of monochromatic complete subgraphs as $(\frac{1}{2}+o(1))\log n$. Adapting previous terminology by which  a 2-colored complete graph of size $n$ with no monochromatic complete subgraph of size $c \log n$ is referred to as $c$-Ramsey,
we say that a 2-colored complete graph $G$ with $n$ vertices is Half-Ramsey if it does not
contain a monochromatic complete subgraph of size $(\frac{1}{2}+o(1))\log n$.

It is not known whether Half-Ramsey graphs exist at all. Here we assume that such graphs do exist, and then study what the property of being Half-Ramsey implies about questions concerning Ramsey multiplicities. One may hope that these implications will either help in actually exhibiting Half-Ramsey graphs, or result in a contradiction that will show that there are no Half-Ramsey graphs.

In Section~\ref{szek_section1} we review known upper bounds on the number of monochromatic complete subgraphs of a Half-Ramsey graph. These upper bounds were derived by Szekely~\cite{DBLP:journals/combinatorica/Szekely84a}. One such bound is the following.

\BTHM
Let $G$ be a Half-Ramsey graph on $n$ vertices. Then $G$ has at most $$n^{\frac{1}{8}(4 - \log e + o(1))\log n} \simeq n^{0.32 \log n}$$ monochromatic complete subgraphs.
\ETHM

In Section~\ref{sec:Szekely2} we provide a matching lower bound, up to low order terms.

\BTHM
Let $G$ be a Half-Ramsey graph on $n$ vertices. Then $G$ has at least $$n^{\frac{1}{8}(4 - \log e - o(1))\log n} \simeq n^{0.32 \log n}$$ monochromatic complete subgraphs.
\ETHM

Moreover, for Half-Ramsey graphs we do not only determine the total number of complete subgraphs, but also determine what we refer to as the {\em profile of Ramsey multiplicities}. Namely, for every $0 \le c \le \frac{1}{2}$ we determine (up to low order terms) the number of monochromatic complete subgraphs of size $c \log n$. This number is $n^{(g_1(c) \pm o(1))\log n}$ for some concave monotonically increasing function $g_1(c)$ that is defined in Section~\ref{sec:profile}. See more details in Section~\ref{sec:profile}.

\subsection{Additional results}

We have some results that shed light on the profile of multiplicities of arbitrary graphs (without requiring them to be Half Ramsey).


In Section \ref{central_results1} we prove a bound that holds for all $t \ge 2$.

\BTHM
For all $t\geq 2$, any $2$-coloring of the edges of the complete graph on $2^{2t-3}$ vertices contains at least
\[
\frac{1}{t!} 2^{\binom{t}{2}-2}
\]
monochromatic complete subgraphs of size $t$.
\ETHM

Furthermore we prove in Section \ref{central_results1} the following variant of Theorem \ref{thm:improveSzekely}.
\BTHM\label{cool_ramsey_thm4_conlon_style_main}
For every $\epsilon>0$ there is a constant $c$ such that the following statement holds. \\
Given a natural number $n$ and a natural number $c \leq b \leq 0.7\log n$.
For
every $2$-coloring $G$ of the edges of the complete graph on $n$ vertices, there is a $\frac{3}{7}b < k < b$ (which depends on $G$)
such that $G$ contains at least
\[
 \frac{n^k}{2^{(1+\epsilon)k^2}}
\]
monochromatic complete subgraphs of size $k$.
\ETHM

In Section \ref{central_results2} we provide improved bounds for sufficiently large $t$.

\BTHM
Let $G$ be a $2$-coloring of a complete graph on $2^t$ vertices. Then for any large enough even $t$, graph $G$ contains at least $2^{\frac{1}{7}t^2}$ monochromatic complete subgraph of size $\frac{1}{2}t$.
\ETHM

Furthermore we prove in Section \ref{central_results2} a bound on the number of monochromatic complete subgraphs of size at most $\frac{1}{2} \log n$.

\BTHM
 Let $G$ be a $2$-coloring of a complete graph on $2^t$ vertices. Then for any large enough $t$, graph $G$ contains at least $2^{c t^2 - O(t\log t) }$ monochromatic complete subgraphs of size at most $t/2$, where  $c = \frac{1}{4}\left(\sqrt{6} - \frac{3}{2}\right) > 0.237$.
\ETHM

%
Another question that we address is the average size of a monochromatic complete subgraphs in a $2$-coloring of a complete graph.
We prove the following in Section \ref{average_section1}.

\BTHM
Let $G$ be a $2$-coloring of a complete graph on $2^t$ vertices. Then for any large enough $t$,
the average size of a monochromatic complete subgraph in $G$ is at least \[\left( 1- \sqrt{\frac{1}{2}} -o(1) \right)t > 0.29t\]
\ETHM

The following upper bound can be easily derived by considering random graphs.

\BTHM
For all large enough $t$, there is a $2$-colorings of a complete graph on $2^t$ vertices, in which the average size of a monochromatic complete subgraph
is at most  $(1+o(1))t$.
\ETHM

Recall that we show that the profile of multiplicities of Half Ramsey graphs has the property that there is a large number of monochromatic complete subgraphs of size roughly $\frac{1}{2}\log n$ (in fact, almost all monochromatic complete subgraphs are of this size), but there is no monochromatic complete subgraph of slightly larger size  $(\frac{1}{2}+o(1))\log n$. As noted above, we do not know if Half Ramsey graphs exist, and hence it is natural to ask whether one can exhibit graphs for which the profile of multiplicities has qualitatively similar properties (e.g., a sudden drop to~0 in the number of monochromatic complete subgraphs). This is one of the motivations for Section \ref{parm_relationship1} that discusses relationships between the maximum size, the average size and the total number of monochromatic complete subgraphs. Among other results, we show:

\BTHM
There is a graph $G$ of size $n$ in which the average size $A(G)$ of a monochromatic complete subgraph satisfies $A(G) = \Theta((\log n)^2)$, and there is no monochromatic complete subgraph of size $A(G) + 2$.
\ETHM


\subsection{Some Notes, Definitions, and Background}
In this paper $\log n$ will denote the binary logarithm, while $\ln n$ will denote the natural logarithm. \\
Some of the related and cited work talk about cliques and independent sets,
which are equivalent to a monochromatic complete subgraph in a $2$-colored complete graph, if we let
one color be edges and the other color be non-edges. Therefore, when we discuss
or use their results we sometimes re-word them to coloring terminologies,
without further comment.

Let the Ramsey Number $R(s, t)$ be the minimum size such that all $2$-colored
graphs of this size have either a blue monochromatic complete subgraph of size $s$ or a red monochromatic complete subgraph of size $t$.
Ramsey's Theorem states that there exists a positive integer $R(s, t)$ such that
this holds. We also state the Erd{\"o}s-Szekeres bound \cite{erdios1935combinatorial}.
\[
R(s,t) \leq \binom{s+t-2}{s-1}
\]
 Let $R(s)$ be the diagonal Ramsey Number, when $s = t$. Simple
known bounds for the diagonal $R(s)$ are of the form:
\[
2^{\frac{s}{2}} \leq R(s) \leq 4^{s-1}
\]

There exist improvements on these bounds in lower order terms, but we do not
use them in this paper so we do not include them here.

Let $k_t(G) $ be the number of monochromatic complete subgraphs of size $t$ in a graph $G$ of
size $n$. Let
\[
k_t(n) = \min \{ k_t(G) : |G|=n \}
\]
Let $c_t(n) = \frac{k_t(n)}{\binom{n}{t}}$, and lastly let $c_t = \lim_{n\rightarrow\infty} c_t(n)$,
so that $c_t$ gives the minimum
fraction of all subsets of size $t$ that are monochromatic complete subgraphs. This notation
is consistent with the related work.

Let $K_i$ denote a clique of size $i$.

\subsection{Related Work}
\label{sec:relatedwork}

The work most related to this paper is that of
Sz{\'{e}}kely~\cite{DBLP:journals/combinatorica/Szekely84a}. He showed that for large enough $t$, in any $2$-coloring of a complete graph on $2^t$ vertices there are at
least $2^{0.1576t^2}$ monochromatic complete subgraphs, and there is $2$-coloring with at most
$2^{(\frac{1}{2}+o(1))t^2}$ monochromatic complete subgraphs. Our Theorem~\ref{thm:improveSzekely} improves over his lower bound. In the same work, Sz{\'{e}}kely also provides upper bounds on Ramsey multiplicities for Half-Ramsey graphs. We prove lower bounds that match his upper bounds (up to low order terms). See Corollary~\ref{cor:HalfRamsey}.

The study of the multiplicity of monochromatic compete subgraphs was introduced by Erd{\"o}s
in 1962 in \cite{erdos1962number}, where  Erd{\"o}s proves that for all graphs,
\begin{equation}\label{erdos_bound_1962_lower_boung1}
c_t \geq \binom{R(t)}{t}^{-1}
\end{equation}
In the same paper,  Erd{\"o}s proves using the probabilistic method that
\begin{equation}\label{erdos_bound_1962_1}
c_t \leq 2^{1-\binom{t}{2}}
\end{equation}
Erd{\"o}s conjectured that the upper bound in Inequality \ref{erdos_bound_1962_1} is tight, or in
other words, that an Erd{\"o}s-R{\'e}nyi $G(n,\frac{1}{2})$ random graph is the graph with the
smallest number of monochromatic complete subgraphs of every size. In 1959, this conjecture
was proved true for the case $t = 3$ by Goodman \cite{goodman1959sets}.

A survey on Ramsey Multiplicity results was published in 1980 by Burr and
Rosta \cite{burr1980ramsey}, in which they extend Erd{\"o}s's conjecture to the multiplicity of any
subgraph, not just monochromatic complete subgraphs.

The conjecture was later disproved by counterexamples in 1989 by Thomason \cite{thomason1989disproof},
who showed that it does not hold for $t \geq 4$. Subsequently, several others worked on upper bounds for $c_t$ for small $t$.
Soon after Thomason's work Franek and R{\"o}dl \cite{franek19932} also gave some different counterexamples based on Cayley
graphs for $t = 4$. Then in 1994,
Jagger, {\v{S}}t'ov{\'\i}{\v{c}}ek and Thomason \cite{jagger1996multiplicities} studied for
which subgraphs the Burr-Rosta conjecture holds, and found that it does not
hold for any graph with $K_4$ as a subgraph, which is consistent with the $t \geq 4$ result found by Thomason.

On the flip side, with regards to the lower bound in Inequality \ref{erdos_bound_1962_lower_boung1}, in 1979
Giraud \cite{giraud1979probleme} proved that $c_4 > \frac{1}{46}$ . More recently in 2012, Conlon~\cite{DBLP:journals/combinatorica/Conlon12}
proved that there must exist at least $\frac{n^t}{C^{(1+o(1))t^2}}$ monochromatic complete subgraphs of size $t$,
in any $2$-colouring of the edges of $K_n$, where $C \approx 2.18$ and $t$ is a constant independent of $n$.
This result is incomparable with our Theorem \ref{cool_ramsey_thm4_conlon_style_main}.

We also give bounds in this paper on the average size of a monochromatic complete subgraph in a graph $G$.
Furthermore we study the ratio between the average size of a monochromatic complete subgraph and the size of a maximum monochromatic complete subgraph in a graph $G$. There
appears to be no previous work directly on this topic.
However, one of the
primary motivations for our work on this topic is from the study of the minimum
of the maximum independent set size over all $K_r$-free graphs of size $n$.

In 1995, Shearer \cite{DBLP:journals/rsa/Shearer95} used the probabilistic method to prove that
$\alpha \geq c'(r)  \frac{n}{d} \frac{\log d}{\log \log (d+1)}$, where $\alpha$ is the size of the maximum independent set and $d$
is the average degree in the graph.
Following his technique, Alon \cite{DBLP:journals/rsa/Alon96} proved
that for a graph in which the neighborhood of every vertex is $r$-colorable, $\alpha \geq \frac{c}{\log(r+1)} \frac{n}{d} \log d$ for some constant $c$. Note
that an $r$-colorable graph is $K_{r+1}$-free, since a clique can contain at most one vertex of each color.

The latest improvement for $K_r$-free graphs is due to Bansal, Gupta and
Guruganesh \cite{DBLP:conf/stoc/BansalGG15}, proving that $\alpha \geq \frac{n}{d} \cdot \max\{\frac{\log d}{r \log \log d},\left( \frac{\log d}{\log r}  \right)^{\frac{1}{2}}\}$.
 There is still a gap in this question, the upper bound being $\frac{n}{d} \frac{\log d}{\log r}$ for $K_r$-free graphs
 (also given in Bansal et al \cite{DBLP:conf/stoc/BansalGG15}). All three of these papers actually prove that a random
independent set in $G$ is of the given size, and then conclude that therefore the
maximum independent set must be at least that size as well.

Thus, knowing the relationship between a random independent set and a
maximum one could be useful in improving these bounds. Analogously we study in this paper the relationship between the average size of a monochromatic complete subgraph and the size of a maximum monochromatic complete subgraph.

\section{The construction of Ramsey Trees}

Let $V=\{v_1,v_2,\ldots,v_n\}$ be the set of vertices of a complete graph $G$ on $n$ vertices. Let {\em red}
and {\em blue} be the two colors of the edges of $G$.
Let $N_{\text{red}}(v)$ be the set of vertices
connected to vertex $v$ by a red edge, and $N_{\text{blue}}(v)$ be the set of vertices connected
to $v$ by a blue edge.
We shall describe several variants of a data structure which we call a Ramsey Tree.
For the sake of clarity, we will use the term \textit{vertex} to
refer to vertices of a graph $G$, and the term \textit{node} to refer to nodes of the respective Ramsey tree.
Estimates on the number of nodes in various levels of the Ramsey tree will allow us to obtain bounds on the number of complete subgraphs of $G$ (see Lemma~\ref{anne_idea1_average_section} for example). The most general variant is the General Ramsey Tree (GRT). Other variants, the Biased Ramsey Tree (BRT) and the restricted Ramsey Tree (RRT), are subtrees of the GRT. They are introduced because their structure is more regular than that of the GRT, and this simplifies the derivation of various estimates that are used in our proofs.

\subsection{Construction of a General Ramsey Tree}

We shall describe how to build a General Ramsey Tree (GRT) $F$ from the graph $G$.
With each node $t$ in the tree, we will associate a vertex $v(t)$
in the graph $G$, the level $l(t)$ of the node, and a bag $B(t)$, where a bag is simply a set of vertices that will
be explained in a few lines. There will be many nodes in the tree construction
associated with a given vertex $v$.

We will build $n$ trees $T_1,T_2,\ldots,T_n$ and later will connect them into one tree $F$. Each tree $T_r$ is rooted at node $r$, where $v(r) = v_r$ so that we have one tree per vertex in $G$.
Each root node $r$ has bag $B(r) = V \backslash v(r)$. Furthermore we set the level of each root node $r$ to be $l(r)=0$.

Each tree $T_r$ is built recursively, as follows.
There is one child of node $t$ for every vertex in the bag $B(t)$, and for each such child $w$ we set the level $l(w) = l(t)+1$.
The
children of node $t$ are split into left and right children.
The left children correspond to vertices in the set $L(t)$ and the right children correspond to vertices in the set $R(t)$,
where $L(t)$ and $R(t)$ are sets of vertices satisfying $|L(t)| + |R(t)| = |B(t)|$ which we define in the following manner:
\begin{itemize}
  \item $L(t)$ contains all the vertices in $B(t)$ that are connected to $v(t)$ in $G$ by a red edge, or in other words:
  $L(t) = N_{\text{red}}(v(t)) \cap B(t)$. For each left child $w$ corresponding to vertex $v(w) \in L(t)$, we let its bag $B(w) = L(t) \backslash v(w)$.
  \item  $R(t)$ contains all the vertices in $B(t)$ that are connected to $v(t)$ in $G$ by a blue edge, or in other words:
  $R(t) = N_{\text{blue}}(v(t)) \cap B(t)$. For each right child $w$ corresponding to vertex $v(w) \in R(t)$, we let its bag $B(w) = R(t) \backslash v(w)$.
\end{itemize}
We apply this recursively, beginning at the root of the tree, then the new
children nodes of the root, then their children, et cetera. The recursions end
when the bags are all empty, and then we have our completed tree $T_r$.

We add a dummy node $d$ (which we denote as the super-root) at level $-1$ and its bag will contain all the vertices of graph $G$ (that is the bag is of size $n$).
The General Ramsey Tree (GRT) $F$ is obtained by connecting the super-root $d$ to the roots of the $n$ trees.

Now we shall describe some properties of the GRT $F$.
Let $Q_i$ be the set of nodes on level $i$ of the GRT $F$. The following two lemmas provide lower and upper bounds on the number of monochromatic complete subgraphs in $G$, as a function of the sizes of $Q_l$ (for various levels $l$).

\BL\label{anne_idea1_average_section}
If for some $l\geq 0$ we have $|Q_l| \geq m$, then graph $G$ contains at least $\sqrt{\frac{m}{(l+1)!}}$ monochromatic complete subgraphs.
\EL
\BPF
The set of vertices corresponding to the nodes in a given path
starting from a root node and ending in level $l$ of the GRT $F$ can appear in at most (l+1)! orders. Hence we have at least $\frac{m}{(l+1)!}$ such paths
where no two paths correspond to the same set of vertices in $G$.
Now consider one of these paths. Every edge in the path is either "going left" or "going right".
In other words, each edge either limits the new bag to a red-connected neighborhood
or a blue-connected neighborhood of the parent node’s corresponding
vertex in $G$. Hence each path induces a red monochromatic complete subgraph in $G$ and a blue monochromatic complete subgraph in $G$ in the following manner:
 We can take one monochromatic complete subgraph to be vertices in $G$ corresponding to the parent nodes of "left going" edges and the second monochromatic complete subgraph consist
 of vertices in $G$ corresponding to the parent nodes of "right going" edges (if the last node in the path has no children we add its corresponding vertex arbitrarily to one of the monochromatic complete subgraphs).
Now notice that no two such paths can induce the same two monochromatic complete subgraphs $C_1$ and $C_2$, as every path corresponds to a different set of vertices in $G$.
Hence we have at least $\sqrt{\frac{m}{(l+1)!}}$ monochromatic complete subgraphs in $G$ and we are done.
\EPF

\BL\label{anne_idea1_average_section2}
If for some $l\geq 0$ we have $|Q_l| \leq m$, then graph $G$ contains at most $\frac{m}{(l+1)!}$ monochromatic complete subgraphs of size $l+1$.
\EL
\BPF
This follows from the fact that every permutation over the vertices of a monochromatic complete subgraph of size $l+1$ appears as a path starting at a root node and ending at level $l$ in the GRT $F$.
\EPF

\BL\label{coollayerslemma1}
For all $i \geq 1$, if $|Q_i|>0$ then $\frac{|Q_{i+1}|}{|Q_i|} \geq \frac{1}{2} \frac{|Q_{i}|}{|Q_{i-1}|} - 1$.
\EL
\BPF
Notice that if $|Q_i|>0$ then $|Q_j|>0$ for all $0 \leq j \leq i$, as each node has a parent in the GRT $F$.
Let $u_1,u_2, \ldots, u_k$ be the nodes on level $i-1$ of the GRT $F$, thus we have $|Q_{i-1}|=k$. Furthermore by the definition of the GRT $F$ we have
\begin{equation}\label{midlayer1}
|Q_i| = \sum_{j=1}^{k} |B(u_j)|
\end{equation}
and
\begin{align}\label{lastlayereq1}
|Q_{i+1}| &= \sum_{j=1}^{k} |L(u_j)|(|L(u_j)|-1) + |R(u_j)|(|R(u_j)|-1) \notag \\
          &\geq \sum_{j=1}^{k} |B(u_j)|\left(\frac{|B(u_j)|}{2}-1\right)
\end{align}
where Inequality (\ref{lastlayereq1}) follows from the fact that $|L(u_j)|(|L(u_j)|-1) + |R(u_j)|(|R(u_j)|-1)$ is minimized when $|L(u_j)| = |R(u_j| = \frac{1}{2}|B(u_j)|$.

Hence
\begin{align}\label{cauchyschwartzargument1}
\frac{|Q_{i+1}|}{|Q_i|} &\geq \frac{\sum_{j=1}^{k} |B(u_j)|\left(\frac{|B(u_j)|}{2}-1\right)}{\sum_{j=1}^{k} |B(u_j)|}
= \frac{1}{2}\frac{\sum_{j=1}^{k} {|B(u_j)|}^2}{\sum_{j=1}^{k} |B(u_j)|} - 1 \notag \\
&\geq \frac{1}{2}\frac{\sum_{j=1}^{k} |B(u_j)|}{k} - 1 \\
&= \frac{1}{2} \frac{|Q_i|}{|Q_{i-1}|} -1 \notag
\end{align}
where Inequality (\ref{cauchyschwartzargument1}) follows from the Cauchy-Schwartz inequality. Thus we are done.
\EPF
\BL\label{important_layer_bound_lem1}
For all $i \geq 0$, if $|Q_i|>0$ then $ \frac{|Q_{i+1}|}{ |Q_i|} > \frac{n}{2^i} - 2 $.
\EL
\BPF
if $|Q_i|>0$ then $|Q_j|>0$ for all $0 \leq j \leq i$, as each node has a parent in the GRT.
Now we will prove the lemma by induction on $i$.
Notice that $|Q_0|=n$ as we have $n$ root nodes in the GRT $F$ and furthermore $|Q_1|=n(n-1)$ as a bag associated with a root node is of size $n-1$.
Hence the base case $i=0$ of the induction follows as $\frac{|Q_1|}{|Q_0|} = n-1 > n-2$. \\
Now assume that the lemma holds for $i-1$, we will prove that the lemma holds for $i$.
By Lemma \ref{coollayerslemma1} we have
\[
\frac{|Q_{i+1}|}{ |Q_i|} \geq \frac{1}{2} \frac{|Q_{i}|}{|Q_{i-1}|} - 1
\]
and by the induction hypothesis we have
\[
\frac{|Q_{i}|}{|Q_{i-1}|} > \frac{n}{2^{i-1}} - 2
\]
Hence we have
\[
\frac{|Q_{i+1}|}{ |Q_i|} > \frac{1}{2}\left(\frac{n}{2^{i-1}} - 2 \right) -1 = \frac{n}{2^i} - 2
\]
and we are done.
\EPF
\BL\label{important_layer_bound_lem2}
If for some $i \geq 0$ and $\delta \geq 0$ , we have $\frac{|Q_{i+1}|}{ |Q_i|} \geq \frac{n^{1+\delta}}{2^{i+1}} - 2 $ then
$\frac{|Q_{j+1}|}{ |Q_j|} \geq \frac{n^{1+\delta}}{2^{j+1}} - 2 $ for all $j \geq i$ for which $|Q_j|>0$.
\EL
\BPF
The proof is almost identical to Lemma \ref{important_layer_bound_lem1} and thus omitted.
\EPF

\subsection{Biased Ramsey Trees}
We start by describing how to build a Biased Ramsey Tree (BRT) $F$ from the graph $G$.

In this construction each node $t$ of the BRT will have a bias parameter $b(t)$ such that $0 \leq b(t) \leq 1$.
Furthermore we always assume that all nodes on the same level of the tree have the same bias (the definition of a level is given in the definition of the BRT below).

As before for GRT, with each node $t$ in the tree, we will associate a vertex $v(t)$
in the graph $G$, the level $l(t)$ of the node, a bag $B(t)$, and the color $c(t)$ of the node (which was implicit in GRTs). In addition, we also associate a bias $0\leq b(t)\leq 1$ with the node, and a parameter $q(t)$ related to $b(t)$.

We will build $n$ trees $T_1,T_2,\ldots,T_n$ (and later connected them into one tree). Each Tree $T_r$ is rooted at node $r$, where $v(r) = v_r$ so that we have one tree per vertex in $G$.
Each root node $r$ has bag  $B(r)$ which is defined as follows: If $|N_{\text{red}}(v_r)| \geq b(r)(n-1)$ then we set $q(r)=b(r)$ and we set $B(r)$ to be an arbitrary subset of $N_{\text{red}}(v_r)$ of size $\lceil q(r)(n-1) \rceil$,
that is $B(r) \subseteq N_{\text{red}}(v_r)$ and $|B(r)| = \lceil q(r)(n-1) \rceil$, furthermore we set the color of the node $c(r)$ to be red. Otherwise we set $q(r)=1-b(r)$ and we set $B(r)$ to be an arbitrary subset of $N_{\text{blue}}(v_r)$ of size $\lceil q(r)(n-1) \rceil$, that is $B(r) \subseteq N_{\text{blue}}(v_r)$ and $|B(r)| = \lceil q(r)(n-1) \rceil$, furthermore we set the color of the node $c(r)$ to be blue. We set the level of each root node $r$ to be $l(r)=0$.

Now we will explain how to build the trees. Each tree is built recursively, as follows.
There is one child of node $t$ for every vertex in the bag $B(t)$ and for each such child $w$ we set the level $l(w) = l(t)+1$.
We define for each child $w$ of node $t$ sets $L(w)$ and $R(w)$ in the following manner:
\begin{itemize}
  \item $L(w)$ contains all the vertices in $B(t)$ that are connected to $v(w)$ in $G$ by a red edge, or in other words:
  $L(w) = N_{\text{red}}(v(w)) \cap B(t)$. 
  \item  $R(w)$ contains all the vertices in $B(t)$ that are connected to $v(w)$ in $G$ by a blue edge, or in other words:
  $R(w) = N_{\text{blue}}(v(w)) \cap B(t)$. 
\end{itemize}
Notice that $|L(w)| + |R(w)| = |B(t)|-1$.
If $|L(w)| \geq b(w)(|B(t)|-1)$ then we set $q(w)=b(w)$ and
we set the bag $B(w)$ to be a subset of $L(w)$ of size $\lceil q(w)(|B(t)|-1) \rceil $, that is  $B(w) \subseteq L(w)$ and $|B(w)| = \lceil q(w)(|B(t)|-1) \rceil$, furthermore we set the color of the node $c(w)$
to be red. Otherwise we set $q(w)=1-b(w)$ and we set $B(w)$ to be a subset of $R(w)$ of size $\lceil q(w)(|B(t)|-1) \rceil$, that is $B(w) \subseteq R(w)$ and $|B(w)| = \lceil q(w)(|B(t)|-1) \rceil$, furthermore we set the color of the node $c(w)$ to be blue.

We apply this recursively, beginning at the root of the tree, then the new
children nodes of the root, then their children, et cetera. The recursions end
when the bags are all empty, and then we have our completed tree $T_r$.

We add a dummy node $d$ (which we denote as the super-root) at level $-1$ and its bag will contain all the vertices of graph $G$ (that is the bag is of size $n$).
The Biased Ramsey Tree (BRT) $F$ consists of these $n$ trees, with roots connected to the super-root.  \smallskip

\subsection{Ramsey Trees (of bias $\frac{1}{2}$)}\label{basic_ramsey_forest_section1}

Biased Ramsey trees in which the bias of all the nodes in the tree is $\frac{1}{2}$ will be particularly convenient for us, especially when the number of vertices in the graph $G$ is a power of~2. Hence we shall reserve the term Ramsey Tree (without mentioning the bias explicitly) to refer to a Biased Ramsey Tree $F$ with bias $\frac{1}{2}$, for a graph $G$ whose number of vertices is $n = 2^q$.
The Ramsey Tree contains $q+1$ levels (not including level $-1$), and the bags at level $q-1$ contain one vertex each.
More generally, the bag size of nodes on level $-1 \leq i < q$ in the Ramsey Tree is $2^{q-i-1}$.

We color the nodes in the final level $q$ of the Ramsey Tree in the following manner.
Let $t$ be a node in level $q$ of the Ramsey Tree. Look at the path from level $0$ up to the parent $p$ of node $t$ (this parent is on level $q-1$). If this path contains at least $\frac{q}{2}$ red nodes we
color node $t$ with red, that is we set $c(t)$ to be red, otherwise we shall set $c(t)$ to be blue.
This coloring ensures us the following fact.
\BL\label{lots_and_lots_of_the_same_color1}
For any $2$-colored complete graph $G$ on $2^q$ vertices,
each path $P$ from a root node to a node in level $q$ in the corresponding Ramsey Tree $F$ contains at least $\lceil \frac{q}{2} \rceil +1$ nodes of the same color.
Furthermore the vertices corresponding to the nodes of the same color in $P$ induce a monochromatic complete subgraph in $G$.
\EL
\BPF
A path $P$ from a root node $u_r$ to a node in level $q$ in the Ramsey Tree contains $q+1$ nodes.
Let $u$ be the last node in path $P$ and let $u_p$ be the parent of $u$ in the path $P$.
The path from $u_r$ ro $u_p$ contains $q$ nodes and thus it contains at least $\lceil \frac{q}{2} \rceil$ nodes of the same color,
assume without loss of generality that this color is red. Thus by the definition of the Ramsey Tree node $u$ will also be colored red.
We conclude that the path $P$ contains at least $\lceil \frac{q}{2} \rceil +1$ nodes of the same color.

Now we shall prove that the vertices corresponding to the nodes of the same color in $P$ induce a monochromatic complete subgraph in $G$.
Let $S$ be a set of nodes of the same color in $P$, assume without loss of generality that this color is red.
Then for any node $u_s \in S$ the vertex $v(u_s)$ is connected to all the vertices in $B(u_s)$ by red edges.
Hence the vertices corresponding to the nodes in $S$ induce a red monochromatic complete subgraph in $G$.
\EPF
\BL\label{ramsey_forest_fact1}
For any $2$-colored complete graph on $2^q$ vertices , the corresponding Ramsey Tree $F$ contains exactly $2^{\binom{q+1}{2}}$ paths
from a root node to a node on level $q$.
\EL
\BPF
Follows from the fact that a bag on level $i$ in the Ramsey Tree is of size $2^{q-i-1}$, hence the number of paths from a root node to a node on level $q$ is
\[
\prod_{i=1}^q 2^i = 2^{\binom{q+1}{2}}.
\]
\EPF

\subsection{Restricted Ramsey Trees}
\label{restricted_ramsey_forest_section2}

Now we describe the construction of Restricted Ramsey Tree (RRT) $F'$ from the graph $G$.
A Restricted Ramsey Tree $F'$ is an induced sub-tree of the  Biased Ramsey Tree $F$ associated with graph $G$, which is defined in the following manner.
Let $Q_0$ be the set of nodes on level $0$ of the BRT $F$.
If $Q_0$ contains more red nodes than blue nodes then level $0$ of $F'$ will contain all the red nodes of $Q_0$,
otherwise level $0$ of $F'$ will contain all the blue nodes of $Q_0$.
Now let $Q_1$ be all the nodes in level $1$ of $F$ which are children (in $F$) of nodes in level $0$ of $F'$.
If $Q_1$ contains more red nodes than blue nodes then level $1$ of $F'$ will contain all the red nodes of $Q_1$,
otherwise level $1$ of $F'$ will contain all the blue nodes of $Q_1$.
And we continue recursively:
let $Q_i$ be all the nodes in level $i$ of $F$ which are children (in $F$) of nodes in level $i-1$ of $F'$.
If $Q_i$ contains more red nodes than blue nodes then level $i$ of $F'$ will contain all the red nodes of $Q_i$,
otherwise level $i$ of $F'$ will contain all the blue nodes of $Q_i$.

Denote the set of nodes in level $i$ of $F'$ by $Q'_i$.
Recall that we will always assume that for each $i$ the bias of all the nodes in $Q'_i$ is the same and we will denote it by $b(i)$.
Hence we have for all $t_1,t_2 \in Q'_i$ that $q(t_1) = q(t_2)$. In particular we can set this value for level $i$ as $q(i)$.
Furthermore by our construction of the RRT we have for each $i$ that the nodes in $Q'_i$ have the same color and we will denote this color by $c(i)$.
Finally we denote by $s(i)$ the size of the bags of the nodes in $Q'_i$ (all such nodes have the same bag size by the construction of the RRT).
\BL\label{tiny_lemma_bag_sizes9}
Given a RRT $F'$, if there is a $p<1$ such that for all $i\geq 0$ we have $q(i) \leq p$ then we have for all $i$ the following.
\[
s(i) \geq n \prod_{j=0}^{i} q(j) - \frac{p}{1-p}
\]
\EL
\BPF
We will prove by induction on $i$.
The bast case $i=0$ follows from the fact that $s(0) \geq q(0)(n-1) \geq q(0) n - p$.
Now we will assume that the claim holds for $i-1$, and will prove it for $i$.
By the definition of the RRT we have
\begin{align}
s(i) &\geq q(i)(s(i-1) -1) \notag \\
&\geq q(i) \left(n \prod_{j=0}^{i-1} q(j) - \frac{p}{1-p} -1\right) &&\text{by the induction hypothesis} \notag \\
&\geq n \prod_{j=0}^{i} q(j) - \frac{p}{1-p} \notag
\end{align}
And thus we are done.
\EPF
Suppose that $l$ is the last level of the RRT $F'$ (that is the bags of nodes on level $l$ are empty) and let $T =\{0,1,\ldots, l\}$.
Recall that $s(-1)=n$ as the bag size of the super-root is $n$.
\BL\label{tiny_lemma_bag_sizes10}
Let $S \subseteq T $ be a set of level indices, where the nodes in all the levels in $S$ of the RRT $F'$ are of the same color,
that is for all $i,j \in S$ we have $c(i)=c(j)$. Then the graph $G$ associated with the RRT $F'$ contains at least
\[
\frac{1}{2^{l+1} (l+1)!} \prod_{i \in S} s(i-1)
\]
monochromatic complete subgraphs of size $|S|$.
\EL
\BPF
Assume $|S|=t$.
Fix a monochromatic complete subgraph $C$ of size $t$.
We shall denote a path from a root node to a node on level $l$ as a full path.
There are at most
\begin{equation}
t! \prod_{i \in T \backslash S} s(i-1)
\end{equation}
different full paths in $F'$ which induce the monochromatic complete subgraph $C$ on the levels of $S$, as there are $t!$ orders in which $C$ can appear on the levels of $S$.
Now by the definition of the RRT the total number of full paths in $F'$ is at least
\begin{align}
\frac{1}{2^{l+1}}\prod_{i \in T} s(i-1)
\end{align}
We conclude that the number of monochromatic complete subgraphs in $G$ of size $t$ is at least
\[
\frac{\prod_{i \in T} s(i-1)}{ 2^{l+1} t! \prod_{i \in T \backslash S} s(i-1) } \geq \frac{1}{2^{l+1} (l+1)!} \prod_{i \in S} s(i-1)
\]
and the proof follows.
\EPF

\section{Monochromatic Complete Subgraphs}\label{central_results1}
One formulation of Ramsey's theorem is the following
\BTHM\label{ramsy_conlon_version}
Any $2$-coloring of the edges of the complete graph on $2^{2t-3}$ vertices contains a monochromatic complete subgraph
of size $t$.
\ETHM
We will start by proving the following strengthening of Theorem \ref{ramsy_conlon_version}.
\BTHM\label{cool_ramsey_theorem1}
For all $t\geq 2$, any $2$-coloring of the edges of the complete graph on $2^{2t-3}$ vertices contains at least
\[
\frac{1}{t!} 2^{\binom{t}{2}-2}
\]
monochromatic complete subgraphs of size $t$.
\ETHM
\BPF
Let $G$ be a $2$-coloring of the edges of the complete graph on $2^{2t-3}$ vertices. Let $F$ be the Ramsey Tree of $G$ as defined in Section \ref{basic_ramsey_forest_section1}.
Notice that by Lemma \ref{lots_and_lots_of_the_same_color1} any path of length $2t-2$ (that is a path from a root node to a node in level $2t-3$)
in $F$ contains either $t$ red nodes or $t$ blue nodes which correspond to
a red or a blue monochromatic complete subgraph of size $t$ in $G$.
Henceforth we will denote a path of length $2t-2$ in $F$ as a full path.

Let $Q=\{0,1,2,\ldots,2t-3\}$. 
Let $S \subseteq Q$ where $S = \{s_1,s_2,\ldots,s_t\}$ and $R = Q \backslash S$ (notice that $|S|=t$ and $|R|=t-2$).
Given a monochromatic complete subgraph $C$ of size $t$ in $G$, it can appear in at most
\begin{equation}\label{basic_bound_fixed_clique1}
t! \prod_{i \in R} 2^i = t! 2^{ \sum_{i \in R} i}
\end{equation}
full paths in $F$ in which the vertices of $C$ correspond to nodes at levels \[\{2t-3-s_1,2t-3-s_2,\ldots,2t-3-s_t\}\]
of the path.
This follows since the sizes of the bags in which the vertices of $C$ appear are simply \[\{2^{s_1},2^{s_2},\ldots,2^{s_t}\}\]
and the product of the sizes of the remaining bags in a full path is $\prod_{i \in R} 2^i$. The vertices of $C$ can appears in $t!$ different orders and thus Equation \ref{basic_bound_fixed_clique1} follows. \\
We conclude from Equation \ref{basic_bound_fixed_clique1} that the number of full paths which contain nodes corresponding to a fixed monochromatic complete subgraph  $C$ of size $t$ is at most
\begin{equation}\label{basic_bound_fixed_clique2}
t! \sum_{ \substack{ R \subseteq Q \\ |R|=t-2}} 2^{ \sum_{i \in R} i}
\end{equation}
The total number of full paths in $F$ is \[\prod_{i=1}^{2t-3} 2^i = 2^{\binom{2t-2}{2}}\]
And hence the number of monochromatic complete subgraphs of size $t$ in $G$ is at least
\begin{align}\label{basic_bound_fixed_clique3}
 2^{\binom{2t-2}{2}}   (t! \sum_{ \substack{ R \subseteq Q \\ |R|=t-2}} 2^{ \sum_{i \in R} i} )^{-1} &=
 2^{\binom{2t-2}{2}}   (t! \sum_{ \substack{ S \subseteq Q \\ |S|=t}} 2^{ \binom{2t-2}{2} -  \sum_{i \in S} i} )^{-1} \\
 &=  (t! \sum_{ \substack{ S \subseteq Q \\ |S|=t}} 2^{-  \sum_{i \in S} i} )^{-1} \notag \\
 &\geq  (t! \sum_{ \substack{ S \subset \mathbb Z_{\ge 0} \\ |S|=t}} 2^{-  \sum_{i \in S} i} )^{-1} \notag
\end{align}
One can prove by induction (see Lemma \ref{tiny_cool_bound1} of Appendix \ref{first_appendix}) that
\begin{equation}\label{basic_bound_fixed_clique4}
\sum_{ \substack{ S \subset \mathbb Z_{\ge 0} \\ |S|=t}} 2^{-  \sum_{i \in S} i} = \frac{2^t}{\prod_{i=1}^{t} {(2^i-1)}  }
\end{equation}
We conclude from Equations \ref{basic_bound_fixed_clique3} and \ref{basic_bound_fixed_clique4}
that the number of monochromatic complete subgraphs of size $t$ in $G$ is at least
\begin{align}\label{basic_bound_fixed_clique5}
\Big(t! \frac{2^t}{\prod_{i=1}^{t} {(2^i-1)}} \Big)^{-1} &= \frac{ \prod_{i=1}^{t} {(2^i-1)}  }{t! 2^t} \\
&= \frac{ \prod_{i=1}^{t} 2^i }{ t! 2^t } \prod_{i=1}^{t} \frac{2^i -1}{2^i} \notag \\
&= \frac{2^{\binom{t}{2}}}{t!} \prod_{i=1}^{t} \Big(1-\frac{1}{2^i}\Big) \notag  \\
&\geq \frac{2^{\binom{t}{2}}}{t!} \prod_{i=1}^{\infty} \Big(1-\frac{1}{2^i}\Big) \notag \\
&\geq \frac{2^{\binom{t}{2}}}{t!} \frac{1}{4} \notag
\end{align}
Where the last inequality follows from the pentagonal number theorem (see Lemma \ref{tiny_cool_bound2} in  Appendix  \ref{first_appendix}).
And thus we are done.
\EPF

In \cite{DBLP:journals/combinatorica/Conlon12} (Theorem $2$) the following is proven.
\BTHM
Let $k,l \geq 1$ be natural numbers. Then, in any red/blue colouring of the edges of $K_n$, there are at least
\[
2^{ -k(l-2) - \binom{k+1}{2} } \binom{n}{k} - O_{k,l}(n^{k-1})
\]
red monochromatic complete subgraphs of size $k$  or at least
\[
2^{ -l(k-2) - \binom{l+1}{2} } \binom{n}{l} - O_{k,l}(n^{l-1})
\]
blue monochromatic complete subgraphs of size $l$.
\ETHM
For $k=l$ Theorem $2$ in \cite{DBLP:journals/combinatorica/Conlon12} implies the following bound.
\BTHM\label{symmetric_conlon_weak}
Let $k \geq 1$ be a natural number. Then, in any red/blue colouring of the edges of $K_n$, there are at least
\[
2^{ \frac{-3k^2 + 3k}{2} } \binom{n}{k} - O_{k}(n^{k-1})
\]
monochromatic complete subgraphs of size $k$.
\ETHM
We prove the following non-asymptotic version of Theorem \ref{symmetric_conlon_weak}.
\BTHM\label{cool_ramsey_cor1}
Given a natural number $n$ and a natural number $2 \leq k \leq \frac{1}{2}\log n$,
Any $2$-coloring of the edges of the complete graph on $n$ vertices contains at least
\[
 2^{\frac{-3k^2 + 5k - 4}{2}} \binom{n}{k}
\]
monochromatic complete subgraphs of size $k$.
\ETHM
\BPF
The proof method is taken from the paper \cite{erdos1962number}.
Let $G$ be a graph on $n$ vertices and suppose $k \leq \frac{1}{2}\log n$.
Let $t=2^{2k-3}$. By Theorem \ref{cool_ramsey_theorem1} every induced subgraph of $t$ vertices in $G$ contains at least
$\frac{1}{k!} 2^{\binom{k}{2}-2}$ monochromatic complete subgraphs of size $k$.
Notice that graph $G$ has $\binom{n}{t}$ induced subgraphs of size $t$. Furthermore each monochromatic complete subgraph of size $k$ is contained in at most
$\binom{n-k}{t-k}$ induced subgraphs of size $t$.
Hence graph $G$ contains at least
\begin{equation}\label{erdos_bound_simple1}
  \frac{\binom{n}{t}}{\binom{n-k}{t-k}} \frac{1}{k!} 2^{\binom{k}{2}-2}
\end{equation}
monochromatic complete subgraphs of size $k$.
Now as $\frac{\binom{n}{t}}{\binom{n-k}{t-k}} = \frac{\binom{n}{k}}{\binom{t}{k}} $ we conclude that graph $G$
contains at least
\begin{align}\label{our_erdos_related_bound_simple2}
  \frac{\binom{n}{t}}{\binom{n-k}{t-k}} \frac{1}{k!} 2^{\binom{k}{2}-2} &= \frac{\binom{n}{k}}{\binom{t}{k}} \frac{1}{k!} 2^{\binom{k}{2}-2}
  \notag \\ &\geq \frac{\binom{n}{k}}{t^k} 2^{\binom{k}{2}-2} \notag \\
  &= \binom{n}{k} 2^{\frac{-3k^2 + 5k - 4}{2}} \notag
\end{align}
monochromatic complete subgraphs of size $k$. And we are done.
\EPF

\BTHM\label{simple_29logn_thm1}
Any $2$-coloring of the edges of the complete graph on $2^t$ vertices contains at least
$
2^{t^2/4 - t\log t -O(t)}
$
monochromatic complete subgraphs of size $r$ where
\[
0.29t < t\left( 1- \sqrt{\frac{1}{2}}\right) \leq r \leq t \sqrt{\frac{1}{2}} < 0.71t
\]
\ETHM
\BPF
Let $G$ be a $2$-coloring of the edges of the complete graph on $2^t$ vertices. Let $F$ be the Ramsey Tree of $G$ as defined in Section \ref{basic_ramsey_forest_section1}.
Henceforth we will denote a path of length $t+1$ (that is a path from a root node to a node in level $t$) in $F$ as a full path.
Set
\[
q= \left\lfloor  t \sqrt{\frac{1}{2}}  \right\rfloor \text{ \quad and \quad }m = \frac{ 2^{\binom{q}{2} - 1 }} { (4t)^t }
\]
Notice that $m= 2^{t^2/4 - t\log t -O(t)}$.
Recall that any set of red nodes in a full path induces a red monochromatic complete subgraph in $G$ and
any set of blue nodes in a full path induces a blue monochromatic complete subgraph in $G$.
If there are at least $m$ different monochromatic complete subgraphs of size $q$ induced by the paths of the Ramsey Tree $F$ we are done.

Thus we may assume that the number of different monochromatic complete subgraphs of size $q$ induced by the paths of the Ramsey Tree $F$ is at most $m$.
Given a fixed monochromatic complete subgraph $C$ of size $q$, there are at most
\begin{equation}\label{ramsey_range_bound1}
b = 2^{\binom{t+1}{2} - \binom{q}{2}} \binom{t+1}{q} q! \leq 2^{\binom{t+1}{2} - \binom{q}{2}} (4t)^t
\end{equation}
different full paths of the Ramsey Tree $F$ which induce the monochromatic complete subgraph $C$. This follows from the fact that we can choose the levels in which nodes corresponding to $C$ appears in $F$ in $\binom{t+1}{q}$ ways and in these levels $C$ can appear in at most $q!$ different orders,
now the product of the sizes of the bags in which the vertices of $C$ do not appear is at most
\[
\prod_{i=q}^{t} 2^i = 2^{\binom{t+1}{2} - \binom{q}{2}}
\]
and bound \ref{ramsey_range_bound1} follows.
By Lemma \ref{ramsey_forest_fact1} the number of full paths in $F$ is $ 2^{\binom{t+1}{2}}$.
Hence if we will remove all the full paths which induce a monochromatic complete subgraph of size $q$ then we will remain with at least
\begin{equation}\label{ramsey_range_bound2}
a= 2^{\binom{t+1}{2}} - m b \geq 2^{\binom{t+1}{2}-1}
\end{equation}
full paths where each such full path induces two monochromatic complete subgraphs $C_1$ and $C_2$
such that
$|C_1|+|C_2| = t+1$ and
$t\left( 1- \sqrt{\frac{1}{2}}\right) \leq |C_1| \leq |C_2|  \leq t \sqrt{\frac{1}{2}}$. \\
The vertices corresponding to a full path can appear in at most (t+1)! orders.
Hence we have at least $\frac{a}{(t+1)!}$ full paths where no two such paths correspond to the same set of vertices in $G$.
Thus no two such paths can induce the same two monochromatic complete subgraphs $C_1$ and $C_2$.  Hence we have at least
\begin{equation}\label{ramsey_square_root1}
\sqrt{\frac{a}{(t+1)!}} = 2^{t^2/4 - t\log t -O(t)}
\end{equation}
monochromatic complete subgraphs
where each such subgraph  is of size at least
$t\left( 1- \sqrt{\frac{1}{2}}\right)$ and at most $t \sqrt{\frac{1}{2}}$.
This concludes the proof.
\EPF
Now we shall show that the bounds in Theorem \ref{simple_29logn_thm1} can be slightly improved (with a more complicated proof).
\BTHM\label{w3logn_thm1}
For large enought $t$, any $2$-coloring of the edges of the complete graph on $2^{2t}$ vertices contains at least
$
2^{(1-o(1))t^2}
$
monochromatic complete subgraphs of size $r$ where
\[
0.6 t \leq r \leq 1.4 t
\]
\ETHM
\BPF
A warning to the reader: it is required in the proof that certain numbers are integers though we have omitted the
corresponding floor/ceiling brackets. Since we do not expect any confusion this will hopefully improve the readability of the proof and of course it will not effect the asymptotics.  \\
Let $G$ be a $2$-coloring of the edges of the complete graph on $n = 2^{2t}$ vertices. Let $F'$ be the Restricted Ramsey Tree (RRT) of $G$ as defined in Section
\ref{restricted_ramsey_forest_section2}, where each node in $F'$ has bias $\frac{1}{2}$.
By the definition of the RRT $F'$ we have that the bags of the nodes in level $i$ of the RRT $F'$ (where $0 \leq i < 2t$) are of size $2^{2t-i-1}$,
that is $s(i) = 2^{2t-i-1}$ (recall also that we defined $s(-1)=2^{2t}$).
Denote by $M_t(G)$ the number of monochromatic complete subgraphs of size $t$ in $G$.
\begin{itemize}
  \item Let $S = \{0,1,\ldots,2t\}$ be the set of the indices of the levels of $F'$.
  \item Let $S_1 = \{0,1,\ldots,t-1\}$ be the set of the indices of the first $t$ levels of $F'$.
  \item Let $S_2 = \{t,t+1,\ldots,2t\}$ be the set of the indices of the last $t+1$ levels of $F'$.
  \item Let $S_R \subseteq S$ be the subset of $S$ of indices of red levels in $F'$, that is for all $i \in S_R$ we have $c(i) = \text{red}$.
  \item Let $S_B \subseteq S$ be the subset of $S$ of indices of  blue levels in $F'$, that is for all $i \in S_B$ we have $c(i) = \text{blue}$.
\end{itemize}
We shall need two more definitions.
\begin{itemize}
\item Given a set $Q \subseteq S$ we define $W(Q) = \prod_{i \in Q} s(i-1)$.
\item Given a set $Q \subseteq S$ we define $W_1(Q) = \log W(Q) = \sum_{i \in Q} \log s(i-1)$.
\end{itemize}
Notice that  $W_1(S) = \sum_{i=0}^{2t} (2t-i) \geq 2 t^2$ and $W_1(S_1) = \sum_{i=0}^{t-1} (2t-i) \geq \frac{3}{2} t^2$. \\
Set $c=0.95$. \\
\textbf{Case $1$:} Assume that $W_1(S_R \bigcap S_1) \leq c t^2$ and $W_1(S_B \bigcap S_1) \leq c t^2$. \\
As $S = S_R \bigcup S_B$ we have that $W_1(S) = W_1(S_R) + W_1(S_B)$ and thus either $W_1(S_R) \geq t^2$ or $W_1(S_B) \geq t^2$. Assume without loss of generality that $W_1(S_R) \geq t^2$.
We will prove in this case that there is a subset $Q \subseteq S_R$ such that $W_1(Q) \geq (1-o(1))t^2$ and $0.6 t \leq |Q| \leq 1.4 t $.
This will be sufficient to prove our theorem in this case as by Lemma \ref{tiny_lemma_bag_sizes10} we will have that $M_{|Q|}(G) \geq 2^{(1-o(1))t^2}$. \\
As $W_1(S_B \bigcap S_1) \leq c t^2$ we have that $W_1(S_R \bigcap S_1) \geq  \left(\frac{3}{2} - c\right)t^2$ (recall that $W_1(S_1) \geq \frac{3}{2} t^2$).
\begin{itemize}
  \item Let $Q_1 \subseteq S_R \bigcap S_1$ be a set of minimum cardinality such that $W_1(Q_1) \geq \left(\frac{3}{2} - c\right)t^2$.
  \item Let $Q_2 \subseteq S_R \setminus Q_1$ be a set of minimum cardinality such that $W_1(Q_2) \geq \left(c - \frac{1}{2}\right)t^2 - 2t$.
  Such a set exists as $W_1(Q_1) \leq \left(\frac{3}{2} - c\right)t^2 + 2t$ by the minimality of $Q_1$ and the fact for all $i$, $\log s(i) \leq 2t$.
  \item Let $Q = Q_1 \bigcup Q_2$. Notice that $W_1(Q) = W_1(Q_1) + W_1(Q_2) \geq (1-o(1))t^2$.
\end{itemize}
Set $t_1 = 0.45 t$ and recall that $c=0.95$. We have
\begin{align}
W_1(\{t-t_1,t-t_1+1,\ldots, t-1 \}) &\geq \sum_{i=t-t_1}^{t-1} (2t-i)
> \frac{1}{2}t_1(t_1+2t) \notag\\
 &> \left(\frac{3}{2} - c\right)t^2 \notag
\end{align}
And thus $|Q_1| \leq t_1$.
Set $t_2 = 0.95 t$. We have
\begin{align}
W_1(\{2t-t_2,2t-t_2+1,\ldots, 2t-1 \}) &\geq \sum_{i=2t-t_2}^{2t-1} (2t-i) >  \frac{1}{2}t_2^2 \notag\\
&> \left(c - \frac{1}{2}\right)t^2 \notag
\end{align}
And thus $|Q_2| \leq t_2$.
We conclude that $|Q| = |Q_1| + |Q_2| \leq 1.4t$. Now we shall prove a lower bound on the cardinality of $Q$. \\
Set $t_3 = 0.55 t$. We have
\begin{align}\label{w30_basic11}
W_1(\{0,1,\ldots,t_3-1\}) &= \sum_{i=0}^{t_3-1} (2t-i)
\leq t_3\left( 2t - \frac{t_3}{2} + 1 \right) \\
 &<  (c-\epsilon)t^2 \notag &&\text{for large enough $t$ and $\epsilon=10^{-9}$}
\end{align}
Set $t_4 = 0.05t$. We have
\begin{align}
W_1(\{t,t+1,\ldots,t+t_4-1\}) &= \sum_{i=t}^{t+t_4-1} (2t-i)
\leq t_4\left( t - \frac{t_4}{2} + 1 \right) \notag\\
 &<  (1-c - \epsilon) t^2 \notag &&\text{for large enough $t$ and $\epsilon=10^{-9}$}
\end{align}
Now notice that for large enough $t$ we have $|Q| \geq t_3 + t_4 = 0.6 t$.
This follows from the following facts:
\begin{enumerate}
  \item Let $Q_3 = Q \bigcap S_1$ and let $Q_4 = Q \bigcap S_2$. Notice that $W_1(Q_3) \leq c t^2$ and thus  $W_1(Q_4) \geq (1-c-o(1)) t^2$.
  \item The cardinality $|Q|=|Q_3|+|Q_4|$ is minimized when  $W_1(Q_3) = (c-o(1)) t^2$ and $W_1(Q_4) = (1-c-o(1)) t^2$.
\end{enumerate}
We have shown that $0.6 t \leq |Q| \leq 1.4 t $ and $W_1(Q) \geq (1-o(1))t^2$. By Lemma \ref{tiny_lemma_bag_sizes10}
we conclude that $M_{|Q|}(G) \geq 2^{(1-o(1))t^2}$ and thus this case is finished. \\
The case we have to consider next is that $W_1(S_R \bigcap S_1) > c t^2$ or $W_1(S_B \bigcap S_1) > c t^2$.
Without loss of generality we will assume that $W_1(S_R \bigcap S_1) > c t^2$ holds. \\
\textbf{Case $2$:} Assume that $W_1(S_R \bigcap S_1) > c t^2$. \\
Set $c'=1.03$ and recall that $c=0.95$. \\
First assume that $W_1(S_R \bigcap S_1) \geq c' t^2$.
Set $t_5 = 0.6 t$. Now as
\begin{align}
W_1(\{0,1,\ldots,t_5-1\}) &= \sum_{i=0}^{t_5-1} (2t-i)
\leq t_5\left( 2t - \frac{t_5}{2} + 1 \right) \notag\\
 &<  (c'-\epsilon)t^2 \notag &&\text{for large enough $t$ and $\epsilon=10^{-9}$}
\end{align}
We conclude that the set $Q = S_R \bigcap S_1$ satisfies $0.6 t = t_5 \leq |Q| \leq t < 1.4 t$ and $W_1(Q) \geq t^2$ and we are done by Lemma \ref{tiny_lemma_bag_sizes10}. \\
For the remainder of the proof we assume that  $c t^2 < W_1(S_R \bigcap S_1) < c' t^2$. \\
We shall build an RRT $F''$ which corresponds to graph $G$ in the following manner.
Levels $0$ up to $t-1$ of the RRT $F''$ will be identical to levels $0$ up to $t-1$ of RRT $F'$.
But in all the remaining levels the bias of the nodes will be $\frac{1}{8}$, that is for RRT $F''$ we have $b(i)=\frac{1}{8}$ for all $i\geq t$.
Let $l$ be the last level of the RRT $F''$. Notice that by Lemma \ref{tiny_lemma_bag_sizes9} we have $l \geq t+\log_8(2^t) - O(1) = \frac{4}{3}t -O(1)$.
Let $l' = \min(2t,l)$.
\begin{itemize}
  \item Let $S' = \{0,1,\ldots,l'\}$  be the set of the indices of the first $l'+1$ levels of $F''$.
  \item Let $S'_1 = \{0,1,\ldots,t-1\}$  be the set of the indices of the first $t$ levels of $F''$.
  \item Let $S'_2 = \{t,t+1,\ldots,l'\}$  be the set of the indices of levels $t$ up to $l'$ of $F''$.
  \item Let $S'_R \subseteq S'$ be the subset of $S'$ of indices of red levels in $F''$, that is for all $i \in S'_R$ we have $c(i) = \text{red}$ in $F''$.
  \item Let $S'_B \subseteq S'$ be the subset of $S'$ of indices of blue levels in $F''$, that is for all $i \in S'_B$ we have $c(i) = \text{blue}$ in $F''$.
  \item Let $Q^1_R = S'_R \bigcap S'_1$ and $Q^2_R = S'_R \bigcap S'_2$.
  \item Let $Q^1_B = S'_B \bigcap S'_1$ and $Q^2_B = S'_B \bigcap S'_2$.
\end{itemize}
We shall need the following definitions. Let $s(i)$ be the bag size in the $i$-th level of the RRT $F''$.
\begin{itemize}
\item Given a set $Q \subseteq S'$ we define $W(Q) = \prod_{i \in Q} s(i-1)$.
\item Given a set $Q \subseteq S'$ we define $W_1(Q) = \log W(Q) = \sum_{i \in Q} \log s(i-1)$.
\end{itemize}
We set $k=0.1 t$ and consider two subcases. \\
\textbf{Subcase $2$-$1$:} Assume that $|Q^2_R| < k$. \\
First we notice that this implies that the index $l$ of the last level of $F''$ satisfies
$l \geq 2t$. This holds as by Lemma \ref{tiny_lemma_bag_sizes9} we have that $s(2t)$ (which is the size of the bags in level $2t$ of the RRT $F''$) satisfies the following bound.
\begin{align}\label{crude_bound_on_level_size7}
s(2t) &\geq n \prod_{j=0}^{2t} q(j) - 7 \notag \\
&\geq 2^{2t} \left(\frac{1}{2}\right)^t \left(\frac{1}{8}\right)^{k}  \left(\frac{7}{8}\right)^{t-k+1} - 7 &&\text{as $|Q^2_R| < k$} \notag \\
&> 2^{t/2} &&\text{for large enough $t$}
\end{align}
 Furthermore bound \ref{crude_bound_on_level_size7} implies that
for $0 \leq i \leq 2t$  the following holds for the bag sizes of RRT $F''$.
\begin{equation}\label{tedious_argument_but_we_need_it7}
s(i) \geq \frac{n}{2} \prod_{j=0}^{i} q(j)
\end{equation}
Let $Z_1$ be an arbitrary subset of $Q^2_B$ of size $t-k=0.9 t$, that is $Z_1 \subseteq Q^2_B$ and $|Z_1|=0.9 t$. \\
Let $Z_2 = Q^1_B \bigcup Z_1$.
As $W_1(S'_R \bigcap S'_1) < c' t^2$ we have that $W_1(Q^1_B) \geq (\frac{3}{2}-c') t^2$.
We also have the following bound.
\begin{align}\label{3log_case1}
  W(Z_1) &=  \prod_{i \in Z_1} s(i-1) \\
  &\geq 2^{-O(t)}  \prod_{i=t+k}^{2t} s(i-1) \notag \\
  &\geq 2^{-O(t)}  \prod_{i=0}^{t-k} 2^{2t} \left(\frac{1}{2}\right)^t \left(\frac{1}{8}\right)^k \left(\frac{7}{8}\right)^i \notag
\end{align}
And this implies that
\begin{align}\label{3log_case1W1}
  W_1(Z_1) &\geq -O(t) + \sum_{i=0}^{t-k} t + \log \left(\frac{1}{8}\right) k + \log \left(\frac{7}{8}\right) i \\
  &\geq -O(t) + (t-k)(t-3k) +  \log \left(\frac{7}{8}\right) \frac{(t-k)^2}{2} \notag \\
  &> 0.55 t^2 &&\text{for large enough $t$} \notag
\end{align}
Hence
\begin{equation}\label{3log_sizebound1}
W_1(Z_2) = W_1(Q^1_B) + W_1(Z_1) \geq (\frac{3}{2}-c')t^2 + 0.55t^2 \geq t^2
\end{equation}
Now we shall prove that $0.6t \leq |Z_2| \leq 1.4t $.
First notice that  \[ |Z_2| \geq |Z_1| = 0.9t > 0.6t \] Furthermore
by Inequality \ref{w30_basic11} we have  $|Q^1_R| \geq t_3 = 0.55t$ and hence $|Q^1_B| \leq t - t_3 = 0.45t$.
Thus \[ |Z_2| = |Q^1_B| + |Z_1| \leq 0.45t + 0.9t < 1.4t \]
We have shown that $0.6t \leq |Z_2| \leq 1.4t $ and that $W_1(Z_2) \geq t^2$. Thus by Lemma \ref{tiny_lemma_bag_sizes10} we have
$M_{|Z_2|}(G) \geq 2^{(1-o(1))t^2}$ and we are done. \\
\textbf{Subcase $2$-$2$:} Assume that $|Q^2_R| \geq k$. \\
Notice that by Lemma \ref{tiny_lemma_bag_sizes10} we have that  for $0 \leq i \leq l-O(1)$  the following holds for the bag sizes of RRT $F''$.
\begin{equation}\label{tedious_argument_but_we_need_it25}
s(i) \geq \frac{n}{2} \prod_{j=0}^{i} q(j)
\end{equation}
Let $Z_3$ be an arbitrary subset of $Q^2_R$ of size $k=0.1 t$, that is $Z_3 \subseteq Q^2_R$ and $|Z_3|=0.1 t$. \\
Let $Z_4 = Q^1_R \bigcup Z_3$.
As $W_1(S'_R \bigcap S'_1) \geq c t^2$ we have that $W_1(Q^1_R) \geq c t^2$.
We also have the following bound.
\begin{align}\label{3log_case3}
  W(Z_3) &=  \prod_{i \in Z_3} s(i-1) \\
  &\geq 2^{-O(t)}  \prod_{i=l'-k}^{l'} s(i-1) \notag \\
  &\geq 2^{-O(t)}  \prod_{i=0}^{k} 2^{2t} \left(\frac{1}{2}\right)^t \left(\frac{7}{8}\right)^{t-k} \left(\frac{1}{8}\right)^i \notag
\end{align}
And this implies that
\begin{align}\label{3log_case1W2}
  W_1(Z_3) &\geq -O(t) + \sum_{i=0}^{k} t + \log \left(\frac{7}{8}\right) (t-k) + \log \left(\frac{1}{8}\right) i\\
  &\geq -O(t) + k t + k(t-k) \log \left(\frac{7}{8}\right) + \log\left(\frac{1}{8}\right) \frac{k^2}{2} \notag \\
  &> 0.06 t^2 &&\text{for large enough $t$} \notag \\
  &> (1-c)t^2 &&\text{as $c=0.95$} \notag
\end{align}
Hence
\begin{equation}\label{3log_sizebound91}
W_1(Z_4) = W_1(Q^1_R) + W_1(Z_3) \geq  c t^2 + (1-c)t^2 = t^2
\end{equation}
Now we shall prove that $0.6t \leq |Z_4| \leq 1.4t $.
First notice that\[ |Z_4| \leq |Q^1_R|+|Z_3| \leq t + 0.1t = 1.1t \]
Thus we have $|Z_4| < 1.4t$. Furthermore
by Inequality \ref{w30_basic11} we have $|Q^1_R| \geq t_3 = 0.55t$.
Thus \[ |Z_4| = |Q^1_R| + |Z_3| \geq 0.55t + 0.1t > 0.6t \]
We have shown that $0.6t \leq |Z_4| \leq 1.4t $ and that $W_1(Z_4) \geq t^2$. Thus by Lemma \ref{tiny_lemma_bag_sizes10} we have
$M_{|Z_4|}(G) \geq 2^{(1-o(1))t^2}$ and we are done.
This concludes the proof.
\EPF
The proof of the following theorem is almost identical to the proof of Theorem \ref{w3logn_thm1} and thus omitted.
\BTHM\label{w3logn_thm2}
Let $\epsilon=10^{-9}$. For large enought $t$, any $2$-coloring of the edges of the complete graph on $2^{2t}$ vertices contains at least
$
2^{(1-o(1))t^2}
$
monochromatic complete subgraphs of size $r$ where
\[
(0.6+\epsilon) t \leq r \leq (1.4-\epsilon) t
\]
\ETHM
\BCR\label{w3logn_cor3}
Any $2$-coloring of the edges of the complete graph on $2^t$ vertices contains at least
$
2^{\left(\frac{1}{4}-o(1)\right)t^2}
$
monochromatic complete subgraphs of size $r$ where
\[
0.3t < r < 0.7t
\]
\ECR
\BCR\label{w3logn_cor4}
Any $2$-coloring of the edges of the complete graph on $n=2^t$ vertices contains at least
\[
 \frac{n^r}{2^{(1+o(1))r^2}}
\]
monochromatic complete subgraphs of size $r$, for some $r$ which satisfies
\[
0.3t < r < 0.7t
\]
\ECR
\BPF
This follows from Corollary \ref{w3logn_cor3}, as we have for $0.3t < r < 0.7t$ that
\[
 \frac{n^r}{2^{(1+o(1))r^2}} \leq 2^{\left(\frac{1}{4}-o(1)\right)t^2}
\]
since the function $f(r)=\frac{n^r}{2^{r^2}}$ attains its maximum at $r=\frac{1}{2} \log n$ and $f(\frac{1}{2}\log n)=2^{t^2/4}$.
\EPF
\BTHM\label{cool_ramsey_thm4_conlon_style}
For every $\epsilon>0$ there is a constant $c$ such that the following statement holds. \\
Given a natural number $n$ and a natural number $c \leq b \leq 0.7\log n$.
For
every $2$-coloring $G$ of the edges of the complete graph on $n$ vertices, there is a $\frac{3}{7}b < k < b$ (which depends on $G$)
such that $G$ contains at least
\[
 \frac{n^k}{2^{(1+\epsilon)k^2}}
\]
monochromatic complete subgraphs of size $k$.
\ETHM
\BPF
We will prove the following equivalent statement.
For every $\epsilon>0$ there is a constant $c$ such that the following statement holds. 
Given a natural number $n$ and a natural number $c \leq t \leq \log n$.
For
every $2$-coloring $G$ of the edges of the complete graph on $n$ vertices, there is a $0.3t < k < 0.7t$ (which depends on $G$)
such that $G$ contains at least
\[
 \frac{n^k}{2^{(1+\epsilon)k^2}}
\]
monochromatic complete subgraphs of size $k$.

Let $G$ be a $2$-coloring of the edges of the complete graph on $n$ vertices and let $t \leq \log n$.
Let $q=2^t$. By Corollary \ref{w3logn_cor3} every induced subgraph of $q$ vertices in $G$ contains at least
$2^{\left(\frac{1}{4}-o(1)\right)t^2}$ monochromatic complete subgraphs of size $0.3t < k < 0.7t$.
Notice that graph $G$ has $\binom{n}{q}$ induced subgraphs of size $q$.
Hence we have for some fixed $0.3t < k < 0.7t$ at least $\frac{1}{t}\binom{n}{q}$ induced subgraphs of size $q$,
all of which contain $2^{\left(\frac{1}{4}-o(1)\right)t^2}$ monochromatic complete subgraphs of size $k$.
Furthermore each monochromatic complete subgraph of size $k$ is contained in at most
$\binom{n-k}{q-k}$ induced subgraphs of size $q$.
Hence graph $G$ contains at least
\begin{equation}\label{erdos_bound_simple1_conlon}
  \frac{1}{t}\frac{\binom{n}{q}}{\binom{n-k}{q-k}} 2^{\left(\frac{1}{4}-o(1)\right)t^2} = \frac{\binom{n}{q}}{\binom{n-k}{q-k}} 2^{\left(\frac{1}{4}-o(1)\right)t^2}
\end{equation}
monochromatic complete subgraphs of size $k$.
Now as $\frac{\binom{n}{q}}{\binom{n-k}{q-k}} = \frac{\binom{n}{k}}{\binom{q}{k}} $ we conclude that graph $G$
contains at least
\begin{align}\label{our_erdos_related_bound_simple2_conlon}
 \frac{\binom{n}{q}}{\binom{n-k}{q-k}} 2^{\left(\frac{1}{4}-o(1)\right)t^2} &= \frac{\binom{n}{k}}{\binom{q}{k}} 2^{\left(\frac{1}{4}-o(1)\right)t^2} \notag \\
 &\geq \frac{n^k}{2^{tk - \left(\frac{1}{4}-o(1)\right)t^2 }} \notag \\
  &\geq  \frac{n^k}{2^{k^2 + o(1)t^2 }} &&\text{as $k^2 \geq tk - \frac{t^2}{4}$ since $\left(k - \frac{t}{2}\right)^2 \geq 0$} \notag \\
   &\geq  \frac{n^k}{2^{(1+o(1))k^2}} &&\text{as $k > 0.3t$}\notag
\end{align}
monochromatic complete subgraphs of size $k$. And we are done.
\EPF

\section{On Half-Ramsey graphs}
\label{sec:HalfRamsey}

Recall the notion of Half-Ramsey graphs from Section~\ref{sec:introduction}.
In Section~\ref{szek_section1} we present upper bounds on the number of monochromatic complete subgraphs of a Half-Ramsey graph. These upper bounds were derived by Szekely~\cite{DBLP:journals/combinatorica/Szekely84a}, and we present their proofs for completeness. In Section~\ref{sec:Szekely2} we prove lower bounds on the number of monochromatic complete subgraphs of a Half-Ramsey graph. Interestingly, our lower bounds match Szekely's upper bounds up to low order terms. The consequence of this is that we can determine with great accuracy the profile of Ramsey multiplicities for Half-Ramsey graphs. This is shown in Section~\ref{sec:profile}.

\subsection{Sz{\'{e}}kely's Bound}\label{szek_section1}

In this section we will describe certain relations between the maximum size of a monochromatic complete subgraph and the number of monochromatic complete subgraphs in a $2$-coloring of a complete graph, first proven in \cite{DBLP:journals/combinatorica/Szekely84a} (we give proofs for completeness).
The following theorem was proven in \cite{DBLP:journals/combinatorica/Szekely84a}.
\BTHM\label{cool_relationship1}
Let $G$ be a $2$-coloring of a complete graph such that $G$ contains no monochromatic complete subgraph of size $t$. Then the number of monochromatic complete subgraphs
in $G$ of size $k<t$ is at most
\[
\frac{2}{k!} \prod_{r=0}^{k-1} \binom{2t-2-r}{t-1}
\]
\ETHM
\BPF
Let $C$ be a red monochromatic complete subgraph of size $r$. A vertex $v$ in $G \backslash C$ is called good if all the edges between $v$ and the vertices of $C$ are red.
Let $L$ be the set of good vertices in $G \backslash C$. We claim that
\begin{equation}\label{good_vertices_lemma1}
  |L| < \binom{2t-2-r}{t-1}
\end{equation}
This holds by the following argument.
Suppose by contradiction that $|L| \geq \binom{2t-2-r}{t-1}$, then by the
Erd\H{o}s-Szekeres bound the complete subgraph in $G$ induced by $L$ contains either a blue monochromatic complete subgraph of size $t$ or a red
monochromatic complete subgraph of size $t-r$. If $G$ contains a blue monochromatic complete subgraph of size $t$ we get a contradiction.
Hence $L$ contains a red monochromatic complete subgraph $R$ of size $t-r$, but then $R \cup C$ is a  monochromatic complete  subgraph of size $t$ in $G$ and we get
a contradiction once again.
Thus the total number of red monochromatic complete subgraphs is at most
\begin{equation}\label{only_red_bound1}
\frac{1}{k!} \prod_{r=0}^{k-1} \binom{2t-2-r}{t-1}
\end{equation}
as we can build every red monochromatic complete subgraph of size $k$ by adding $k$ vertices iteratively
where at stage $r$ the number of vertices that can be added to the red monochromatic complete subgraph of size $r$ is at most $\binom{2t-2-r}{t-1}$.
Furthermore each such red monochromatic complete subgraph of size $k$ will appear in $k!$ different orders.
The total number of blue monochromatic complete subgraphs is also bounded by \ref{only_red_bound1} and hence the theorem follows.
\EPF
Furthermore the following was proven in \cite{DBLP:journals/combinatorica/Szekely84a} (once again we give a proof for completeness).
\BCR\label{cool_relationship2}
Let $G$ be a $2$-coloring of a complete graph such that $G$ contains no monochromatic complete subgraph of size $t$. Then the number of monochromatic complete subgraphs
in $G$ is at most
\[
2^{(2-\frac{\log e}{2})t^2 + O(t \log t)}
\]
\ECR
\BPF
By Theorem \ref{cool_relationship1} the number of monochromatic complete subgraphs of $G$ of size $k<t$
is bounded from above by
\[
2 \prod_{r=0}^{k-1} \binom{2t-2-r}{t-1}
\]
Hence the total number of monochromatic complete subgraphs of $G$ is bounded from above by
\begin{align}\label{simplifying_a_bit_cool_relationship2}
  2t \prod_{r=0}^{t-1} \binom{2t-2-r}{t-1} &\leq  2t \prod_{r=0}^{t} \binom{2t-r}{t} \\
  &= 2t \prod_{r=0}^{t} \binom{t+r}{t} \notag
\end{align}
Set $G(n) = \prod_{k=0}^{n} k!$ (the function $G(n)$ is known in the literature as the Barnes $G$-Function).
Then we have
\begin{align}\label{asymptotic_stuff1_cool_relationship2}
  2t \prod_{r=0}^{t} \binom{t+r}{t} &\leq 2t \frac{\prod_{r=t}^{2t} r!}{\prod_{r=0}^{t} r! \cdot (t!)^t } \\
  &= 2t \frac{\prod_{r=0}^{2t} r!}{\prod_{r=0}^{t} r! \cdot \prod_{r=0}^{t-1} r! \cdot (t!)^t } \notag \\
  & =2t \left(\frac{1}{t!}\right)^t \frac{G(2t)}{G(t)G(t-1)} \notag \\
  &\leq  2t \left(\frac{e}{t}\right)^{t^2} \frac{G(2t)}{G(t)G(t-1)}  && \text{from the Stirling bound $t! \geq \left(\frac{t}{e}\right)^t $ } \notag
\end{align}
Now by Equation $(A.6)$ of \cite{voros1987spectral} we have the following bound
\begin{equation}\label{voros_striling1}
 e^{n^2\left( \frac{\ln n}{2} -\frac{3}{4}\right) - O(n \ln n)} \leq  G(n) \leq e^{n^2\left( \frac{\ln n}{2} -\frac{3}{4}\right) + O(n \ln n)}
\end{equation}
We conclude from Inequalities \ref{asymptotic_stuff1_cool_relationship2} and \ref{voros_striling1} that
\[
2t \prod_{r=0}^{t} \binom{t+r}{t} \leq 2^{(2-\frac{\log e}{2})t^2 + O(t \log t)}
\]
and thus we are done.
\EPF
\BCR
For large enough $t$, Let $G$ be a $2$-coloring of a complete graph such that $G$ contains no monochromatic complete subgraph of size $0.5001 t$. Then the number of monochromatic complete subgraphs
in $G$ is less than $2^{0.32t^2}$
\ECR
Another way to formulate the corollary above is that for any large enough $t$, if $G$ is a $2$-coloring of a complete graph such that $G$ contains at least  $2^{0.32t^2}$ monochromatic complete subgraphs
then $G$ contains a monochromatic complete subgraph of size at least $0.5001 t$.

Now we shall present Sz{\'{e}}kely's bound in full generality.
Define a function $g(c)$ for real  $0 \leq c \leq 1$ in the following manner:
\[
g(c) =
\begin{cases}
 \frac{1}{2}(4-c\log e + (1-c)^2 \log(1-c) - (2-c)^2 \log(2-c) ) ,          &\text{if $0 \leq c < 1$} \\
 \frac{1}{2}(4-\log e)  &\text{if $c = 1$}
\end{cases}
\]
The function $g(c)$ is monotonically increasing in the domain $[0,1]$ and we have $g(0)=0$ and $g(1)= \frac{1}{2}(4-\log e)$.
Furthermore the function $g(c)$ is concave in the domain $[0,1]$.
\BTHM\label{full_sezekely_relationship1}
Let $G$ be a $2$-coloring of a complete graph such that $G$ contains no monochromatic complete subgraph of size $t$.
Let $0 \leq c \leq 1$ be a constant. Then the number of monochromatic complete subgraphs
in $G$ of size at most $c t$ is at most
\[
2^{ g(c) t^2 + O(t \log t)}
\]
\ETHM
\BPF
We may assume that $c t$ is an integer.
By Theorem \ref{cool_relationship1} the number of monochromatic complete subgraphs of $G$ of size at most $c t$ is bounded from above by
\begin{align}\label{szekfull_simplifying_a_bit_cool_relationship2}
  2t \prod_{r=0}^{c t} \binom{2t-2-r}{t-1} &\leq  2t \prod_{r=0}^{c t} \binom{2t-r}{t}
\end{align}
Set $G(n) = \prod_{k=0}^{n} k!$. Then we have

\begin{align}\label{szekfull_asymptotic_stuff1_cool_relationship2}
 2t \prod_{r=0}^{c t} \binom{2t-r}{t} &= 2t \prod_{r=0}^{c t} \frac{(2t-r)!}{t!(t-r)!} \\
  &\leq 2t \frac{G(2t)}{G((2-c)t-1)} \frac{1}{(t!)^{c t}} \frac{G((1-c)t)}{G(t)} \notag \\
 &\leq 2t \frac{G(2t)}{G((2-c)t-1)} \frac{G((1-c)t)}{G(t)}  \left(\frac{e}{t}\right)^{c t^2}  && \text{from the Stirling bound $t! \geq \left(\frac{t}{e}\right)^t $ } \notag \\
  &\leq 2t \frac{G(2t)}{G((2-c)t)} \frac{G((1-c)t)}{G(t)}  \left(\frac{e}{t}\right)^{c t^2} (2t)^{2t} \notag
\end{align}
We can conclude from Inequalities \ref{szekfull_asymptotic_stuff1_cool_relationship2} and \ref{voros_striling1} that
\[
2t \prod_{r=0}^{c t} \binom{2t-r}{t} \leq 2^{g(c) t^2 + O(t \log t)}
\]
and thus we are done.
\EPF
\BCR\label{szekfullbound11}
For large enough $n$, Let $G$ be a $2$-coloring of a complete graph such that $G$ contains no monochromatic complete subgraph of size $(\frac{1}{2} + o(1))\log n$. Then
for any constant $0 \leq c \leq \frac{1}{2}$
the number of monochromatic complete subgraphs of size at most $c \log n$ in $G$ is at most $n^{\left(\frac{1}{4}g(2c) +o(1)\right) \log n}$.
\ECR
\BPF
The proof follows from setting $t= (\frac{1}{2} + o(1))\log n$ in Theorem \ref{full_sezekely_relationship1}.
\EPF
Notice that the following consequence of Corollary \ref{szekfullbound11} follows trivially.
\BCR
For large enough $n$, Let $G$ be a $2$-coloring of a complete graph such that $G$ contains no monochromatic complete subgraph of size $(\frac{1}{2} + o(1))\log n$. Then
for any constant $0 \leq c \leq \frac{1}{2}$
the number of monochromatic complete subgraphs of size exactly $c \log n$ in $G$ is at most $n^{\left(\frac{1}{4}g(2c) +o(1)\right) \log n}$.
\ECR

\subsection{A converse to Sz{\'{e}}kely's bound}
\label{sec:Szekely2}

In this section we will prove the following theorem.
\BTHM\label{perhaps_the_most_important_theorem}
Let $G$ be a $2$-coloring of the edges of the complete graph on $2^{2t}$ vertices such that the maximum size of a monochromatic complete subgraph in graph $G$ is $(1+o(1))t$.
Then at least one of the following two statements holds:
\begin{enumerate}
  \item $G$ contains at least $2^{\frac{1}{2}(4 - \log e -o(1))t^2}$ red monochromatic complete subgraphs of size at least $t$.
   Furthermore $G$ contains at least $2^{\frac{1}{2}(4 - \log e -o(1))t^2}$ blue monochromatic complete subgraphs of size at least $(1-o(1))t$.
  \item $G$ contains at least $2^{\frac{1}{2}(4 - \log e -o(1))t^2}$ blue monochromatic complete subgraphs of size at least $t$.
   Furthermore $G$ contains at least $2^{\frac{1}{2}(4 - \log e -o(1))t^2}$ red monochromatic complete subgraphs of size at least $(1-o(1))t$.
\end{enumerate}
\ETHM
The remainder of this section is devoted to the proof of the theorem above.
We will need the Erd{\"o}s-Szekeres bound \cite{erdios1935combinatorial} in this section.
\BTHM\label{fundumental_erdos_bound1}
\[
R(s,t) \leq \binom{s+t-2}{s-1}
\]
\ETHM
Furthermore we will requires the following estimates.
\BL\label{ln_estimate1}
If $x>0$ then $\ln(1+x) > x - \frac{1}{2} x^2$.
\EL
\BL\label{ln_estimate2}
If $0< x< 0.69$ then $\ln(1-x) > -x - x^2$.
\EL
\BL\label{entropy_estimate1}
Suppose that $n p$ is an integer in the range $[0, n]$, then $ \binom{n}{n p} \leq 2 ^{n H(p)}$.
Where $H(p) = -p\log(p) - (1-p)\log(1-p)$ is the binary entropy function.
\EL
Henceforth for all of this section $G$ will be a $2$-coloring of the edges of the complete graph on $2^{2t}$ vertices.
We will denote this $2$-coloring as graph $G$. Let $F'$ be the Restricted Ramsey Tree (RRT) of $G$ as defined in Section
\ref{restricted_ramsey_forest_section2}, where each node in $F'$ has bias $\frac{1}{2}$.
By the definition of the RRT $F'$ we have that the bags of the nodes in level $i$ of the RRT $F'$ (where $-1 \leq i < 2t$) are of size $2^{2t-i-1}$,
that is $s(i) = 2^{2t-i-1}$. 
Recall that $c(i)$ is the color of the nodes in level $i$ of the RRT $F'$.
We will use the following notation.
\begin{itemize}
  \item $c_r(i)$ denotes the number of red nodes in a path from a root node to a node in level $i$ of the RRT $F'$ (including the nodes on the root level and on level $i$).
  \item $c_b(i)$ denotes the number of blue nodes in a path from a root node to a node in level $i$ of the RRT $F'$ (including the nodes on the root level and on level $i$).
  \item $c'_r(i)$ denotes the number of red nodes in a path from a node in level $i$ to a node in level $2t$ of the RRT $F'$
  (including the nodes on level $i$ and on level $2t$).
  \item $c'_b(i)$ denotes the number of blue nodes in a path from a node in level $i$ to a node in level $2t$ of the RRT $F'$ (including the nodes on level $i$ and on level $2t$).
\end{itemize}
Notice that $c_r(i) + c_b(i) = i+1$ and $c'_r(i) + c'_b(i) = 2t-i+1$.

Next we will prove Lemma \ref{szekely_bag_size_lem1} , Lemma \ref{balancing_lemma1} and Lemma \ref{balancing_lemma_cor2} which are needed to show
that for each $i$ we have $c_r(i) \simeq c_b(i) \simeq \frac{i}{2}$ and $c'_r(i) \simeq c'_b(i) \simeq \frac{2t - i}{2} $. This statement will be made precise in Lemma \ref{balancing_lemma_cor2}.
\BL\label{szekely_bag_size_lem1}
If $G$ contains no monochromatic complete subgraph of size $q+1$ then for all $i\geq 0$ the bag size in RRT $F'$ satisfies
\[s(i-1) <  \binom{2q-c_b(i-1)-c_r(i-1)}{q-c_b(i-1)} =  \binom{2q-i}{q-c_b(i-1)}\]
\EL
\BPF
Let $u$ be a node in level $i-1$ of the RRT $F'$. Assume by contradiction that  \[ s(i-1) \geq \binom{2q-c_b(i-1)-c_r(i-1)}{q-c_b(i)} \]
Then by Theorem \ref{fundumental_erdos_bound1} the vertices of $B(u)$ contain either a red monochromatic complete subgraph of size $q-c_r(i-1)+1$ or a
blue monochromatic complete subgraph of size  $q-c_b(i-1)+1$.

If $B(u)$ contain a red monochromatic complete subgraph of size $q-c_r(i-1)+1$
then this set together with the vertices associated with the $c_r(i-1)$ red nodes in a path from a root node to node $u$ induce a
red monochromatic complete subgraph of size $q+1$ and we reach a contradiction.

If $B(u)$ contain a blue monochromatic complete subgraph of size $q-c_b(i-1)+1$
then this set together with the vertices associated with the $c_b(i-1)$ blue nodes in a path from a root node to node $u$ induce a
blue monochromatic complete subgraph of size $q+1$ and we reach a contradiction.
Thus we are done.
\EPF
\BL\label{balancing_lemma1}
Assume that the maximum size of a monochromatic complete subgraph in $G$ is $(1+\epsilon)t$ where $0 < \epsilon<10^{-4}$.
Set $d=3 \epsilon^{\frac{1}{4}}$. Let $i = 2\alpha t$ for $d \leq \alpha \leq  1 - d$.
Then $c_b(i-1) \geq (\alpha - d)t$ and $c_r(i-1) \geq (\alpha - d)t$.
\EL
\BPF
Let $q=(1+\epsilon)t$. Recall that $i = 2\alpha t$ and we can assume that $\sqrt{\epsilon} \leq \alpha \leq  1 - \sqrt{\epsilon}$ by the choice of $d$.
 Let $c_b(i-1) = \beta t$. By the definition of the RRT $F'$ we have $s(i-1) = 2^{2t-i}$. On the other hand by Lemma \ref{szekely_bag_size_lem1} we have
 $s(i-1) < \binom{2q-i}{q-c_b(i-1)}$. Thus we have
 \begin{equation}\label{balancing_act1}
 2^{2t-i} < \binom{2q-i}{q-c_b(i-1)}
 \end{equation}
 And we conclude from Equation \ref{balancing_act1} that
 \begin{align}\label{balancing_act2}
 2^{2t(1-\alpha)} &< \binom{2t(1+\epsilon-\alpha)}{t(1+\epsilon-\beta)} \\
 &\leq  2^{2t(1+\epsilon-\alpha) H\left(\frac{1}{2} \frac{1+\epsilon-\beta}{1+\epsilon-\alpha}\right)} &&\text{by Lemma \ref{entropy_estimate1}} \notag
 \end{align}
Where $H(p) = -p\log(p) - (1-p)\log(1-p)$ is the binary entropy function.
Hence from Inequality \ref{balancing_act2}  we conclude that
\begin{align}\label{balancing_act3}
 1 &<    \left(\frac{1+\epsilon-\alpha}{1-\alpha}\right) H\left(\frac{1}{2} \frac{1+\epsilon-\beta}{1+\epsilon-\alpha}\right) \\
   &=  \left(1+\frac{\epsilon}{1-\alpha}\right) H\left(\frac{1}{2} \frac{1+\epsilon-\beta}{1+\epsilon-\alpha}\right) \notag \\
   &\leq  \left(1+\sqrt{\epsilon}\right) H\left(\frac{1}{2} \frac{1+\epsilon-\beta}{1+\epsilon-\alpha}\right) \notag &&\text{as $\alpha \leq 1 - \sqrt{\epsilon}$}
 \end{align}
Now suppose that $\beta = \alpha - \delta$ for some $\delta > 0$ (if $\delta \leq 0$ then we are done).
Hence we have from Inequality \ref{balancing_act3} the following.
\begin{align}\label{balancing_act4}
1 &< \left(1+\sqrt{\epsilon}\right) H\left(\frac{1}{2}\left( 1 + \frac{\delta}{1+\epsilon-\alpha} \right) \right) \\
&\leq \left(1+\sqrt{\epsilon}\right) H\left(\frac{1}{2}\left( 1 + \delta \right) \right) \notag
 \end{align}
Where the last inequality follows from the following two facts.
\begin{enumerate}
  \item As $\alpha \geq \epsilon$ we have $\delta \leq \frac{\delta}{1+\epsilon-\alpha}$.
  \item The function $H\left(\frac{1}{2}\left( 1 + \delta \right) \right)$ is monotonically decreasing for $\delta\geq 0$.
\end{enumerate}
Now multiplying Inequality \ref{balancing_act4} by $(1-\sqrt{\epsilon})$ we get
\begin{equation}\label{balancing_act5}
H\left(\frac{1}{2}\left( 1 + \delta \right) \right) \geq 1 - \sqrt{\epsilon}
\end{equation}
By the definition of the binary entropy function we have
\begin{align}\label{balancing_act6}
H\left(\frac{1}{2}\left( 1 + \delta \right) \right) &= -\frac{1+\delta}{2} \log\left(\frac{1+\delta}{2}\right) -\frac{1-\delta}{2} \log\left(\frac{1-\delta}{2}\right) \\
&= 1 - \frac{1}{2}\log(1-\delta^2) - \frac{\delta}{2} \log\left( \frac{1+\delta}{1-\delta}\right) \notag \\
&=  1 - \frac{1}{2}\log(1-\delta^2) - \frac{\delta}{2} \log\left(1+\frac{2\delta}{1-\delta}\right) \notag \\
&\leq  1 - \frac{1}{2}\log(1-\delta^2) - \frac{\delta}{2} \log\left(1+2\delta\right) \notag
\end{align}
We conclude from Inequalities \ref{balancing_act6} and \ref{balancing_act5} that
\[
 \frac{1}{2}\log(1-\delta^2) + \frac{\delta}{2} \log\left(1+2\delta\right) \leq \sqrt{\epsilon}
\]
And thus in particular
\begin{equation}\label{balancing_act7}
\frac{1}{2}\ln(1-\delta^2) + \frac{\delta}{2} \ln\left(1+2\delta\right) \leq \sqrt{\epsilon}
\end{equation}
Now by Lemma \ref{ln_estimate1} and \ref{ln_estimate2} we have
\begin{align}\label{balancing_act8}
 \frac{1}{2}\ln(1-\delta^2) + \frac{\delta}{2} \ln\left(1+2\delta\right) &\geq
 \frac{1}{2}(-\delta^2 - \delta^4) + \frac{\delta}{2} (2\delta - \frac{1}{2} (2\delta)^2) \\
 &= \frac{1}{2}\delta^2 (1 - \delta^2 - 2\delta) \notag  \\
 &\geq  \frac{1}{8}\delta^2 \notag
\end{align}
Where the last inequality follows from the fact that as $\epsilon < 10^{-4}$ we have from Inequality \ref{balancing_act5} that $\delta < \frac{1}{4}$,
and thus $1 - \delta^2 - 2\delta \geq \frac{1}{4}$.
We conclude from Inequalities \ref{balancing_act8} and \ref{balancing_act7} that
\[
\frac{1}{8}\delta^2 \leq \sqrt{\epsilon}
\]
Hence we have
\[\delta \leq 3 \epsilon^{\frac{1}{4}}\]
And thus  \[ c_b(i-1) = \beta t = (\alpha - \delta)t \geq (\alpha - d)t \]
Furthermore we can show that $c_r(i-1) \geq (\alpha - d)t$ in the exact same way as
 \begin{equation}\label{balancing_act100}\notag
 2^{2t-i} < \binom{2q-i}{q-c_r(i-1)}
 \end{equation}
And thus we are done.
\EPF
\BL\label{balancing_lemma_cor2}
Assume that the maximum size of a monochromatic complete subgraph in $G$ is $(1+\epsilon)t$ where $0 < \epsilon<10^{-4}$.
Set $d_1 = 5 \epsilon^{\frac{1}{4}}$. Then for all $2 d_1 t \leq i \leq (1-d_1)2t$ and large enough $t$ (depending on $\epsilon$ but not on $i$) the following inequalities holds:
\begin{enumerate}
  \item  $\frac{i}{2} - d_1 t \leq c_r(i-1) \leq \frac{i}{2} + d_1 t$
  \item  $\frac{i}{2} - d_1 t \leq c_b(i-1) \leq \frac{i}{2} + d_1 t$
  \item  $t - \frac{i}{2} - d_1 t \leq c'_r(i+1) \leq t - \frac{i}{2} + d_1 t $
  \item  $t - \frac{i}{2} - d_1 t \leq c'_b(i+1) \leq t - \frac{i}{2} + d_1 t $
\end{enumerate}
\EL
\BPF
Set $d=3 \epsilon^{\frac{1}{4}}$. Recall that we assume that $2 d_1 t \leq i \leq (1-d_1)2t$.
We start by proving parts 1 and 2. By Lemma \ref{balancing_lemma1} we have  $\frac{i}{2} - d t \leq c_r(i-1)$ and
$\frac{i}{2} - d t \leq c_b(i-1)$. As $c_b(i-1) + c_r(i-1) = i$ we conclude that $c_r(i-1) \leq \frac{i}{2} + d t$ and $c_b(i-1) \leq \frac{i}{2} + d t$.

Now we shall prove parts 3 and 4. As the maximum size of a monochromatic complete subgraph in $G$ is $(1+\epsilon)t$ we have that  $c_r(2t) \geq (1-\epsilon)t$ and
$c_b(2t) \geq (1-\epsilon)t$. Hence we have
\begin{align}
c'_r(i+1) &\geq c_r(2r) - c_r(i-1) -1 \\ \notag &\geq (1-\epsilon)t - (\frac{i}{2} + d t) -1 \notag \\ &\geq t - \frac{i}{2} - d_1 t && \text{for large enough $t$}
\end{align}
In the same manner we can prove that
$
 c'_b(i+1) \geq t - \frac{i}{2} - d_1 t
$.
Now as $c'_r(i+1) + c'_b(i+1) = 2t-i$ we conclude that $c'_r(i+1) \leq t - \frac{i}{2} + d_1 t$ and $c'_b(i+1) \leq t - \frac{i}{2} + d_1 t $.
And thus we are done.
\EPF

Henceforth we will denote a path starting from a root node and ending in level $2t$  as a full path.
\BD
We will define an auxiliary bipartite graph $H(W,Z)$ in the following manner:
Each vertex $w \in W$ correspond to a set of $c_r(2t)$ vertices associated with the red nodes of some full path in the RRT $F'$.
Each vertex $z \in Z$ correspond to a set of $c_b(2t)$ vertices associated with the blue nodes of some full path in the RRT $F'$.
Finally there is an edge between vertex $w \in W$ and $z \in Z$ if and only if there is a full path in the RRT $F'$ such that the set of  vertices associated with the nodes of this path is the union
of the sets of vertices associated with $w$ and $z$.
\ED
Recall that the set of nodes in level $i$ of the RRT $F'$ is denoted by $Q'_i$.
By the definition of the RRT $F'$ we have that
\begin{equation}\label{old_ineq_once_more1}
|Q'_{2t}| \geq \frac{ \prod_{j=0}^{2t} 2^{2t-j}   }{2^{2t+1}} =  2^{2t^2 - t - 1}
\end{equation}

The vertices corresponding to the nodes in a given full path in $F'$ can appear in at most (2t+1)! orders.
Hence by bound \ref{old_ineq_once_more1}
we have at least
 \[
  \frac{ 2^{2t^2 - t - 1} } { (2t+1)! } =  2^{2t^2 - O(t\log t)}
 \]
full paths with distinct vertices.
Hence we have the following
\BL\label{Hedgenum1}
The number of edges in the auxiliary graph $H$ is at least
\[
2^{2t^2 - O(t\log t)}
\]
\EL
Denote by $S_R$ the set of indices of the levels in the RRT $F'$ where red nodes appear and denote by
$S_B$ the set of indices of the levels in the RRT $F'$ where blue nodes appear.
Notice that $S_R \bigcup S_B = \{0,1,\ldots,2t \}$.
\begin{itemize}
  \item Denote by $d_H(w)$ the degree of a vertex $w \in W$ of graph $H$.
  \item Denote by $d_H(z)$ the degree of a vertex $z \in Z$ of graph $H$.
  \item Denote by $d_H^W$  the maximum degree in graph $H$ of the vertices in  $W$.
  \item Denote by $d_H^Z$ the maximum degree in graph $H$ of the vertices in $Z$.
  \item Denote by $E(H)$ the number of edges in graph $H$
\end{itemize}
By the definition of the auxiliary graph $H$ we have the following bound.
\BL\label{uri_lem1}
The number of red monochromatic complete subgraphs in $G$ of size $c_r(2t)$ is at least
\[ |W| \geq \frac{E(H)}{d_H^W} \]
The number of blue monochromatic complete subgraphs in $G$ of size $c_b(2t)$ is at least
\[ |Z| \geq \frac{E(H)}{d_H^Z} \]
\EL
\BL\label{Hdegree_lem1}
Assume that $G$ contains no monochromatic complete subgraph of size $q+1$.
Let $w \in W$ be a vertex in the auxiliary graph $H(W,Z)$, then we have
\[
d_H(w) < (2t+1)!\prod_{i \in S_B} \binom{2q - c'_b(i+1) - c_r(i-1) } {q - c_r(i-1) }
\]
\EL
\BPF
We will give a bound on the number of full paths of the RRT $F'$ in which vertices which are associated with the red nodes of the path are fixed to be the vertices associated with vertex $w$ of graph $H$ in some order. Notice that there are at most $(t + o(t))! \leq (2t+1)!$ such orders.

Now given $i \in S_B$ assume that all vertices associated with blue nodes on levels $i < j \in S_B$ in the path are fixed.
We claim that the number of different vertices of graph $G$
which can be associated with nodes on level $i \in S_B$ in such path is less than
\[
\binom{2q - c'_b(i+1) - c_r(i-1) } {q - c_r(i-1) }
\]
This holds as each set of $\binom{2q - c'_b(i+1) - c_r(i-1) } {q - c_r(i-1) }$ vertices in graph $G$
contains by Theorem \ref{fundumental_erdos_bound1} either a red monochromatic complete subgraph of size $q -  c_r(i-1) +1 $ or a blue
monochromatic complete subgraph of size $q - c'_b(i+1) +1$. \\
If it contains a red monochromatic complete subgraph of size $q -  c_r(i-1) +1 $ then this set together with the $c_r(i-1)$ vertices associated with the red nodes of index at most $i-1$
induce a red monochromatic complete subgraph of size $q+1$ in $G$ by contradiction to our assumption. \\
If it contains a blue monochromatic complete subgraph of size $q -  c'_b(i+1) +1 $ then this set together with the $ c'_b(i+1)$ vertices associated with the blue nodes of index at least $i+1$
induce a blue monochromatic complete subgraph of size $q+1$ in $G$ by contradiction to our assumption.
We conclude that we can complete the full path in which the red nodes are fixed in at most
\[\prod_{i \in S_B} \binom{2q - c'_b(i+1) - c_r(i-1) } {q - c_r(i-1) } \]
different ways and thus
\[
d_H(w) < (2t+1)!\prod_{i \in S_B} \binom{2q - c'_b(i+1) - c_r(i-1) } {q - c_r(i-1) }
\]
and we are done.
\EPF
\BL\label{Hdegree_lem2}
Assume that $G$ contains no monochromatic complete subgraph of size $q+1$.
Let $z \in Z$ be a vertex in the auxiliary graph $H(W,Z)$, then we have
\[
d_H(z) < (2t+1)!\prod_{i \in S_R} \binom{2q - c'_r(i+1) - c_b(i-1) } {q - c_b(i-1) }
\]
\EL
\BPF
The proof is almost identical to the proof of Lemma \ref{Hdegree_lem1} and thus ommited.
\EPF
We will need the following technical lemmas.
\BL\label{binom_estimate1}
For all non-negative integers $n,k,t$ we have $\binom{n+t}{k+t} \leq \binom{n}{k} (\frac{n+t}{k})^t$.
\EL
\BPF
Notice that for $t\geq 1$ we have $\binom{n+t}{k+t} = \binom{n+t-1}{k+t-1} (\frac{n+t}{k+t})$.
Our Lemma follows by iterating this identity.
\EPF
\BL\label{zigzagftw1}
Assume that the maximum size of a monochromatic complete subgraph in $G$ is $q=(1+\epsilon)t$  where $0 < \epsilon<10^{-4}$.
Set $d_2 = 6 \epsilon^{\frac{1}{4}}$. Then for all
\[ 4 d_2 t \leq i \leq (1-4d_2)t  \text{ and } (1 + 8 d_2) t \leq i \leq (2-8d_2)t \]
and large enough $t$ (depending on $\epsilon$ but not on $i$) the following inequality holds
\[
\binom{2q - c'_r(i+1) - c_b(i-1) } {q - c_b(i-1) } \leq \binom{t}{\lceil \frac{i}{2}\rceil} 2^{- \log(d_2) \lceil 2 d_2 t \rceil}
\]
\EL
\BPF
By the symmetry of the binomial coefficient $ \binom{t}{\lceil \frac{i}{2}\rceil} $ it is sufficient to prove the claim in the range
$  4 d_2 t \leq i \leq (1-4d_2)t $.
Set $d_1 = 5 \epsilon^{\frac{1}{4}}$. By Lemma \ref{balancing_lemma_cor2} we have
\begin{equation}\label{zigzag_calc1}
2q - c'_r(i+1) - c_b(i-1) \leq 2(1+\epsilon)t - (t -\frac{i}{2} - d_1 t) - (\frac{i}{2} - d_1 t)\leq \lceil t + 2 d_2 t \rceil
\end{equation}
Furthermore we have by Lemma \ref{balancing_lemma_cor2} that
\begin{align}\label{zigzag_calc2}
q - c'_r(i+1) &\leq (1+\epsilon)t - (t - \frac{i}{2} - d_1 t) \leq \left \lceil \frac{i}{2} + d_2 t \right \rceil \\
&\leq \left \lceil \frac{t}{2} -d_2 t \right \rceil &&\text{as $i \leq  (1-4d_2)t$} \label{little_binom_help1}
\end{align}
We conclude from Inequalities \ref{zigzag_calc1} and \ref{zigzag_calc2} that
\begin{align}\label{hop}
\binom{2q - c'_r(i+1) - c_b(i-1) } {q - c_b(i-1) } &\leq \binom{\lceil t + 2 d_2 t \rceil}{  \lceil \frac{i}{2} + d_2 t  \rceil } \notag  \\
 &\leq \binom{\lceil t + 2 d_2 t \rceil}{  \lceil \frac{i}{2} + 2 d_2 t  \rceil } &&\text{by Inequality \ref{little_binom_help1}} \notag \\
 &\leq \left(\frac{4t}{i}\right)^{\lceil 2d_2 t \rceil} \binom{t}{\lceil \frac{i}{2}\rceil } \notag &&\text{by Lemma \ref{binom_estimate1}} \notag \\
 &\leq \left(\frac{1}{d_2}\right)^{\lceil 2d_2 t \rceil} \binom{t}{\lceil \frac{i}{2}\rceil }  \notag  &&\text{as $i \geq 4 d_2 t$} \\
 &= \binom{t}{\lceil \frac{i}{2}\rceil} 2^{- \log(d_2) \lceil 2 d_2 t \rceil} \notag
\end{align}
And thus we are done.
\EPF
\BL\label{zigzagftw2}
Assume that the maximum size of a monochromatic complete subgraph in $G$ is $q=(1+\epsilon)t$  where $0 < \epsilon<10^{-4}$.
Set $d_2 = 6 \epsilon^{\frac{1}{4}}$. Then for all
\[ 4 d_2 t \leq i \leq (1-4d_2)t  \text{ and } (1 + 8 d_2) t \leq i \leq (2-8d_2)t \]
and large enough $t$ (depending on $\epsilon$ but not on $i$) the following inequality holds
\[
\binom{2q - c'_b(i+1) - c_r(i-1) } {q - c_r(i-1) } \leq \binom{t}{\lceil \frac{i}{2}\rceil} 2^{- \log(d_2) \lceil 2 d_2 t \rceil}
\]
\EL
\BPF
The proof is almost identical to the proof of Lemmat \ref{zigzagftw1} and thus omitted.
\EPF
Let $G(n) = \prod_{k=0}^{n} k!$.
We will need the following estimate from  \cite{voros1987spectral}.
\BL\label{superfactorial_estimate13}
\begin{equation}
 e^{n^2\left( \frac{\ln n}{2} -\frac{3}{4}\right) - O(n \ln n)} \leq  G(n) \leq e^{n^2\left( \frac{\ln n}{2} -\frac{3}{4}\right) + O(n \ln n)} \notag
\end{equation}
\EL
\BTHM
Assume that the maximum size of a monochromatic complete subgraph in graph $G$ is $(1+o(1))t$ and
recall that graph $G$ contains $2^{2t}$ vertices.
Then at least one of the following two statements holds:
\begin{enumerate}
  \item $G$ contains at least $2^{\frac{1}{2}(4 - \log e -o(1))t^2}$ red monochromatic complete subgraphs of size at least $t$.
   Furthermore $G$ contains at least $2^{\frac{1}{2}(4 - \log e -o(1))t^2}$ blue monochromatic complete subgraphs of size at least $(1-o(1))t$.
  \item $G$ contains at least $2^{\frac{1}{2}(4 - \log e -o(1))t^2}$ blue monochromatic complete subgraphs of size at least $t$.
   Furthermore $G$ contains at least $2^{\frac{1}{2}(4 - \log e -o(1))t^2}$ red monochromatic complete subgraphs of size at least $(1-o(1))t$.
\end{enumerate}
\ETHM
\BPF
By Lemma \ref{Hdegree_lem1} and Lemma \ref{Hdegree_lem2} we have
\begin{align}
d_H^W \cdot d_H^Z &\leq ((2t+1)!)^2 \prod_{i \in S_B} \binom{2q - c'_b(i+1) - c_r(i-1) } {q - c_r(i-1) } \cdot \prod_{i \in S_R} \binom{2q - c'_r(i+1) - c_b(i-1) } {q - c_b(i-1) }
\notag \\
&\leq 2^{ o(t^2) } \prod_{i=0}^{2t} \binom{t}{\lceil \frac{i}{2}\rceil} \text{\qquad \qquad \qquad \qquad \qquad \qquad \quad by Lemma \ref{zigzagftw1} and Lemma \ref{zigzagftw2}} \notag \\
&\leq  2^{ o(t^2) } \left (\prod_{i=0}^{t} \binom{t}{i} \right )^2 \notag \\
&\leq 2^{ o(t^2) } \frac{(t!)^{2t}}{G(t)^4} \notag \\
&\leq e^{(1+o(1))t^2}  \text{\qquad \qquad \qquad \qquad \quad \quad by Lemma \ref{superfactorial_estimate13} and the Stirling approximation} \notag
\end{align}
We conclude that
\begin{equation}\label{productZW_upper_bound1}
d_H^W \cdot d_H^Z \leq  2^{(\log(e)+o(1))t^2}
\end{equation}
Recall that by Lemma \ref{Hedgenum1} the number of edges in the auxiliary graph $H$ satisfies $E(H) \geq 2^{(2-o(1))t^2}$.
By Lemma \ref{uri_lem1} we have that the number of red monochromatic complete subgraphs in $G$ of size $c_r(2t)$ is at least
\[ \frac{E(H)}{d_H^W} \]
And the number of blue monochromatic complete subgraphs in $G$ of size $c_b(2t)$ is at least
\[ \frac{E(H)}{d_H^Z} \]
Hence by Corollary \ref{cool_relationship2} 
we have
\begin{equation}
  \frac{E(H)}{d_H^W} \leq 2^{(2-\frac{1}{2}\log e +o(1))t^2}  \notag
\end{equation}
Now as by Lemma \ref{Hedgenum1} we have $E(H) \geq 2^{(2-o(1))t^2}$ we conclude that
\begin{equation}\label{much_needed_lowerbound1}
  d_H^W \geq 2^{(\frac{1}{2}\log(e)- o(1))t^2}
\end{equation}
And by the same manner we have
\begin{equation}\label{much_needed_lowerbound2}
  d_H^Z \geq 2^{(\frac{1}{2}\log(e)- o(1))t^2}
\end{equation}
Combining Inequality \ref{productZW_upper_bound1} with Inequalities \ref{much_needed_lowerbound1} and \ref{much_needed_lowerbound2} we conclude that
\begin{equation}\label{much_needed_lowerbound3}
  d_H^W \leq 2^{(\frac{1}{2}\log(e)+ o(1))t^2}
\end{equation}
And that
\begin{equation}\label{much_needed_lowerbound4}
  d_H^Z \leq 2^{(\frac{1}{2}\log(e)+ o(1))t^2}
\end{equation}
Thus the number of red monochromatic complete subgraphs of size $c_r(2t)$ is at least
\[
 \frac{E(H)}{d_H^W} \geq 2^{(2  - \frac{1}{2}\log e -o(1))t^2}
\]
And the number of blue monochromatic complete subgraphs of size $c_b(2t)$ is at least
\[
 \frac{E(H)}{d_H^Z} \geq 2^{(2  - \frac{1}{2}\log e -o(1))t^2}
\]
And the proof is finished taking into account that $c_r(2t) + c_b(2t)=2t+1$ and the fact that
$c_b(2t) \leq (1+o(1))t$  and $c_r(2t) \leq (1+o(1))t$.
\EPF

\subsection{The profile of Ramsey multiplicities of Half-Ramsey graph}
\label{sec:profile}

\begin{figure}[h]
  \includegraphics[width=\linewidth]{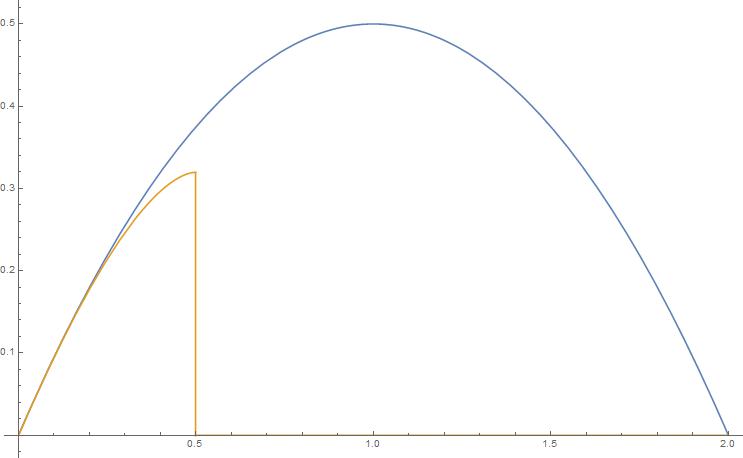}
  \caption{Functions $g_1(c)$ and $g_2(c)$.}
  \label{fig:boat1}
\end{figure}

We say that a graph $G$ with $n$ vertices is Half-Ramsey if it does not
contain either a clique or an independent set of size $(\frac{1}{2}+o(1))\log n$.
This notation is similar to the usual notation in the literature where  a graph
$G$ of size $n$ is denoted as $c$-Ramsey, if it has no clique or independent set of size
$c \log n$.
\BTHM\label{halftheorem1}
Let $G$ be a Half-Ramsey graph on $2^{2t}$ vertices. Then $G$ contains at least $2^{\frac{1}{2}(4 - \log e -o(1))t^2}$ monochromatic complete subgraphs of size at least $t$.
Furthermore $G$ contains at most  $2^{\frac{1}{2}(4 - \log e + o(1))t^2}$ monochromatic complete subgraphs in total. Notice that the upper and lower bounds are identical up to the $o(1)$ term in the exponent.
\ETHM
\BPF
The upper bound follows from Corollary \ref{cool_relationship2}. While the lower bound follows from Theorem \ref{perhaps_the_most_important_theorem}.
\EPF
Now we shall prove a more general result.
Recall that we  defined a function $g(c)$ for real  $0 \leq c \leq 1$ in the following manner in Section \ref{szek_section1}.
\[
g(c) =
\begin{cases}
 \frac{1}{2}(4-c\log e + (1-c)^2 \log(1-c) - (2-c)^2 \log(2-c) ) ,          &\text{if $0 \leq c < 1$} \\
 \frac{1}{2}(4-\log e)  &\text{if $c = 1$}
\end{cases}
\]
Recall that the function $g(c)$ is monotonically increasing in the domain $[0,1]$ and we have $g(0)=0$ and $g(1)= \frac{1}{2}(4-\log e)$.
Furthermore the function $g(c)$ is concave in the domain $[0,1]$.
\BTHM
Let $G$ be a Half-Ramsey graph on $2^{2t}$ vertices.
Let $0 \leq c \leq 1$ be a constant. Then the number of monochromatic complete subgraphs
in $G$ of size $c t$ is at most
\[
2^{ (g(c)+o(1)) t^2}
\]
and at least
\[
2^{ (g(c)-o(1)) t^2}
\]
\ETHM
\BPF
The upper bound follows from Theorem \ref{full_sezekely_relationship1}.
Now we shall prove the lower bound.
Denote by $C_k(G)$ the number of monochromatic complete subgraph of size $k$ in $G$.
Set $k=c t$. First we notice that the total number of monochromatic complete subgraphs in $G$ is at most
\begin{equation}\label{halfeq1}
2 C_k(G) \prod_{r=k}^{t+o(t)} \binom{2t+o(t)-r}{t+o(t)}
\end{equation}
The proof of this bound is almost identical to the proof of Lemma \ref{cool_relationship1} and thus omitted.
Now by Theorem \ref{halftheorem1} the total number of monochromatic complete subgraphs in $G$ is at least
\begin{equation}\label{halfeq2}
2^{\frac{1}{2}(4 - \log e -o(1))t^2}
\end{equation}
Hence from the Bounds \ref{halfeq1} and \ref{halfeq2} we get the following inequality.
\[
2 C_k(G) \prod_{r=k}^{t+o(t)} \binom{2t+o(t)-r}{t+o(t)} \geq 2^{\frac{1}{2}(4 - \log e -o(1))t^2}
\]
Applying the asymptotics of the Barnes $G$-Function as done in Theorem \ref{full_sezekely_relationship1} we get
\[
C_k(G) \geq 2^{ (g(c)-o(1)) t^2}
\]
and thus we are done.
\EPF
Let
\[
g_1(c) =
\begin{cases}
 \frac{1}{4}g(2c)  ,          &\text{if $0 \leq c \leq \frac{1}{2}$} \\
 0  &\text{if $c > \frac{1}{2}$}
\end{cases}
\]
\BCR
\label{cor:HalfRamsey}
Let $G$ be a Half-Ramsey graph on $n$ vertices. Then
for any constant $0 \le c \le \frac{1}{2}$,
the number of monochromatic complete subgraphs of size exactly $c \log n$ in $G$ is at most
\[ n^{\left(g_1(c) +o(1)\right) \log n} \]
and at least
\[ n^{\left(g_1(c) -o(1)\right) \log n} \]
\ECR
\BCR
Let $G$ be a Half-Ramsey graph on $n$ vertices. Then the average size of a monochromatic complete subgraph of $G$ is at least
\[
\left(\frac{1}{2} - o(1)\right)\log n
\]
\ECR
\BPF
This follows from the fact that $g_1(c)$ is monotonically increasing in the range $[0,\frac{1}{2}]$ (see Figure \ref{fig:boat1}).
\EPF
Now we can ask ourselves what is the analogous theorem when $G$ is a random graph $G(n,\frac{1}{2})$ (that is each edge is chosen independently with probability $\frac{1}{2}$).
One can prove the following in this case using standard techniques.
Let $g_2(c) = c - \frac{1}{2}c^2$.
\BTHM
Let $G(n,\frac{1}{2})$ be a random graph in which each edge is chosen with probability $\frac{1}{2}$.
For any constant $0 \leq c < 2$ we have asymptotically almost surely that
the number of monochromatic complete subgraphs of size exactly $c \log n$ in $G(n,\frac{1}{2})$ is at most
\[ n^{\left(g_2(c) +o(1)\right) \log n} \]
and at least
\[ n^{\left(g_2(c) -o(1)\right) \log n} \]
\ETHM
Figure \ref{fig:boat1} shows the functions $g_1(c)$ and $g_2(c)$.

\section{Monochromatic Complete Subgraphs revisited}\label{central_results2}

In this section we prove additional results (beyond those of Section~\ref{central_results1}) concerning the number of monochromatic complete subgraphs of various sizes.

\BTHM
Let $c = \sqrt{6} - \frac{3}{2} > 0.949 $. For large enough $t$,
any $2$-coloring of the edges of the complete graph on $2^{2t}$ vertices contains at least
$
2^{c t^2 - O(t\log t) }
$
different monochromatic complete subgraphs of size at most $t$.
\ETHM
\BPF
Let $G$ be a $2$-coloring of the edges of the complete graph on $2^{2t}$ vertices. Let $F'$ be the Restricted Ramsey Tree (RRT) of $G$ as defined in Section
\ref{restricted_ramsey_forest_section2}, where each node in $F'$ has bias $\frac{1}{2}$.

Recall that the set of nodes in level $i$ of $F'$ is denoted by $Q'_i$.
By the definition of the RRT $F'$ we have that the bags of the nodes in level $i$ of the RRT $F'$ (where $0 \leq i < 2t$) are of size $2^{2t-i-1}$,
that is $s(i) = 2^{2t-i-1}$  (we also assume that $s(-1) = 2^{2t}$).
Hence we have for all $0 \leq i \leq 2t$ the following bound.
\begin{equation}\label{new_ramsey_forest_level_bound1}
|Q'_i| \geq \frac{ \prod_{j=0}^{i} 2^{2t-j}   }{2^{i+1}} = 2^{\frac{(4t-i)(i+1)}{2} -i-1} =  2^{2ti - \frac{1}{2}i^2 - O(t)}
\end{equation}
Recall that $c(i)$ is the color of the nodes in level $i$ of the RRT $F'$.
We define the color order of $F'$ as a vector $D=(c(0),c(1),\ldots,c(2t))$.
Let $a$ the smallest index for which the prefix of $D$, which we denote by $D' = (c(0),c(1),\ldots,c(a))$, contains $t$ entries of the same color.
Notice that $t-1 \leq a \leq 2t-2$. \\
Denote by $M_t(G)$  the number of monochromatic complete subgraphs of size at most $t$ in $G$
By Lemma \ref{tiny_lemma_bag_sizes10} and the fact that in the worst case the $t$ levels of the same color are levels $a,a-1,\ldots,a-(t-1)$
we have that
\begin{align}\label{important_bound_on_less_then_equal_size2}
M_t(G) &\geq 2^{-O(t \log t)} \prod_{i=a-t+1}^{a} s(i-1) \\
&\geq  2^{-O(t \log t)} \prod_{i=a-t+1}^{a} 2^{2t-i} &&\text{(as $s(i) = 2^{2t-i-1}$)} \notag \\
&\geq  2^{ \frac{5}{2}t^2 - at - O(t \log t) } \notag
\end{align}
Now we shall give another lower bound on $M_t(G)$.
Let $F''$ be the induced sub-tree of $F'$ which contains only levels $0$ up to $a$ of $F'$.
Henceforth we will denote a path starting from a root node and ending in level $a$  as a full path.
By bound \ref{new_ramsey_forest_level_bound1} tree $F''$ contains at least
\begin{equation}\label{number_of_full_paths_in_forest7}
  2^{2ta - \frac{1}{2}a^2 - O(t)}
\end{equation} full paths.
Notice that each full path in $F''$  induces two monochromatic complete subgraphs $C_1$ and $C_2$ such that
$|C_1|\leq t$ and $|C_2| \leq t$.
The vertices corresponding to the nodes in a given full path in $F''$ can appear in at most (a+1)! orders.
Hence by bound \ref{number_of_full_paths_in_forest7}  and the fact that $a < 2t$
we have at least
 \[
 \frac{1}{(a+1)!} 2^{2ta - \frac{1}{2}a^2 - O(t)} \geq  2^{2ta - \frac{1}{2}a^2 - O(t \log t)}
 \]
full paths with distinct vertices.
Thus no two such paths can induce the same two monochromatic complete subgraphs $C_1$ and $C_2$.  Hence we have at least
\begin{equation}\label{ramsey_square_root_less_than_equal2}
\sqrt{ 2^{2ta - \frac{1}{2}a^2 - O(t \log t)} }=  2^{ta - \frac{1}{4}a^2 - O(t \log t)}
\end{equation}
monochromatic complete subgraphs
where each such subgraph  is of size at most $t$.
We conclude that
\begin{equation}\label{important_bound_on_less_then_equal_size3}
  M_t(G) \geq   2^{ta - \frac{1}{4}a^2 - O(t \log t)  }
\end{equation}
Notice that bound \ref{important_bound_on_less_then_equal_size3} is monotonically increasing for $0 \leq a \leq 2t$, while
bound \ref{important_bound_on_less_then_equal_size2} is monotonically decreasing for $a\geq 0$. \\
Thus for $a \geq (4-\sqrt{6})t$ we get from bound  \ref{important_bound_on_less_then_equal_size3} that $M_t(G) \geq 2^{(\sqrt{6} - \frac{3}{2})t^2 - O(t\log t) } $.
Furthermore for $a \leq (4-\sqrt{6})t$ we get from bound  \ref{important_bound_on_less_then_equal_size2} that $M_t(G) \geq 2^{(\sqrt{6} - \frac{3}{2})t^2 - O(t\log t) } $.
And thus we are done.
\EPF
\BCR
Let $c = \frac{1}{4}(\sqrt{6} - \frac{3}{2}) > 0.237$. For large enough $t$,
any $2$-coloring of the edges of the complete graph on $2^{t}$ vertices contains at least
$
2^{c t^2 - O(t\log t) }
$
different monochromatic complete subgraphs of size at most $t/2$.
\ECR
\BTHM
For large enough $t$,  any $2$-coloring of the edges of the complete graph on $2^{2t}$ vertices contains at least
$
2^{\frac{4}{7} t^2}
$
monochromatic complete subgraphs of size $t$.
\ETHM
\BPF
Let $G$ be a $2$-coloring of the edges of the complete graph on $n = 2^{2t}$ vertices. Let $F'$ be the Restricted Ramsey Tree (RRT) of $G$ as defined in Section
\ref{restricted_ramsey_forest_section2}, where each node in $F'$ has bias $\frac{1}{2}$.
By the definition of the RRT $F'$ we have that the bags of the nodes in level $i$ of the RRT $F'$ (where $0 \leq i < 2t$) are of size $2^{2t-i-1}$,
that is $s(i) = 2^{2t-i-1}$ (recall also that we defined $s(-1)=2^{2t}$).
Denote by $M_t(G)$ the number of monochromatic complete subgraphs of size $t$ in $G$.
\begin{itemize}
  \item Let $S = \{0,1,\ldots,2t\}$ be the set of the indices of the levels of $F'$.
  \item Let $S_1 = \{0,1,\ldots,t-1\}$ be the set of the indices of the first $t$ levels of $F'$.
  \item Let $S_R \subseteq S$ be the subset of $S$ of indices of red levels in $F'$, that is for all $i \in S_R$ we have $c(i) = \text{red}$.
  \item Let $S_B \subseteq S$ be the subset of $S$ of indices of  blue levels in $F'$, that is for all $i \in S_B$ we have $c(i) = \text{blue}$.
\end{itemize}
\textbf{Case $1$:} Assume that $|S_R \bigcap S_1| \geq \frac{t}{3}$ and $|S_B \bigcap S_1| \geq \frac{t}{3}$. \\
Since $|S| \geq 2t$  we have either $|S_R \bigcap S| \geq t$  or $|S_B \bigcap S| \geq t$.
Assume without loss of generality that $|S_R \bigcap S| \geq t$.
Recall that $|S_R \bigcap S_1| \geq \frac{t}{3}$. Let $Q \subseteq S_R$ be a subset of $S_R$ such that $|Q|=t$ and $|Q \bigcap S_1| \geq \frac{t}{3}$.
Now by Lemma
\ref{tiny_lemma_bag_sizes10} we have
\begin{align}\label{large_monoclique_eqn1}
  M_t(G) &\geq  2^{-O(t \log t)}  \prod_{i \in Q} s(i-1) \\
  &\geq 2^{-O(t \log t)} \prod_{\frac{2t}{3} < i < t} s(i-1) \prod_{\frac{4t}{3} < i < 2t} s(i-1) &&\text{as  $|Q \bigcap S_1| \geq \frac{t}{3}$}   \notag \\
  &\geq 2^{-O(t \log t)} \prod_{\frac{2t}{3} < i < t} 2^{2t-i} \prod_{\frac{4t}{3} < i < 2t} 2^{2t-i} \notag \\
  &\geq  2^{-O(t \log t)} \: 2^{\frac{11}{18}t^2 - O(t) } \notag \\
  &\geq  2^{\frac{4}{7}t^2} &&\text{for large enough $t$} \notag
\end{align}
This concludes Case $1$ of the proof. Hence we may assume without loss of generality that $|S_B \bigcap S_1| < \frac{t}{3}$ and thus $|S_R \bigcap S_1| > \frac{2t}{3}$. \\
\textbf{Case $2$:} Assume that $|S_B \bigcap S_1| < \frac{t}{3}$ and $|S_R \bigcap S_1| > \frac{2t}{3}$. \\
We shall build an RRT $F''$ which corresponds to graph $G$ in the following manner.
Levels $0$ up to $t-1$ of the RRT $F''$ will be identical to levels $0$ up to $t-1$ of RRT $F'$.
But in all the remaining levels the bias of the nodes will be $\frac{2}{5}$, that is for RRT $F''$ we have $b(i)=\frac{2}{5}$ for all $i\geq t$.
Let $l$ be the last level of the RRT $F''$.
\begin{itemize}
  \item Let $S' = \{0,1,\ldots,l\}$ be the set of the indices of the levels of $F''$.
  \item Let $S'_1 = \{0,1,\ldots,t-1\}$ be the set of the indices of the first $t$ levels of $F''$.
  \item Let $S'_2 = \{t,t+1,\ldots,l\}$ be the set of the indices of levels $t$ up to $l$ of $F''$.
  \item Let $S'_R \subseteq S'$ be the subset of $S'$ of indices of red levels in $F''$, that is for all $i \in S'_R$ we have $c(i) = \text{red}$ in $F''$.
  \item Let $S'_B \subseteq S'$ be the subset of $S'$ of indices of blue levels in $F''$, that is for all $i \in S'_B$ we have $c(i) = \text{blue}$ in $F''$.
\end{itemize}
First we assume that $|S'_R| \geq t$. Let $Q' \subseteq S'_R$ be a subset of $S'_R$ such that $|Q'|=t$ and $|Q' \bigcap S'_1| \geq \frac{2t}{3}$.
Thus by Lemma \ref{tiny_lemma_bag_sizes10} (applied for RRT $F''$) we have
\begin{align}
M_t(G) &\geq 2^{-O(t \log t)}  \prod_{i \in Q'} s(i-1) \\
  &\geq  2^{-O(t \log t)}  \prod_{\frac{t}{3} < i < t} s(i-1) &&\text{as $|Q' \bigcap S'_1| \geq \frac{2t}{3}$ }\notag \\
&\geq 2^{-O(t \log t)}  \prod_{\frac{t}{3} < i < t} 2^{2t-i} \notag \\
&\geq  2^{-O(t \log t)}\: 2^{\frac{8}{9}t^2 - O(t) } \notag \\
&\geq  2^{\frac{4}{7}t^2} &&\text{for large enough $t$} \notag
\end{align}
Thus henceforth we may assume that $|S'_R| < t$. This implies that $|S'_R \bigcap S'_2| < \frac{t}{3}$.
First we notice that this implies that the index $l$ of the last level of $F''$ satisfies
$l \geq 2t$. This holds as by Lemma \ref{tiny_lemma_bag_sizes9} we have that $s(2t)$ (which is the size of the bags in level $2t$ of the RRT $F''$) satisfies the following bound.
\begin{align}\label{crude_bound_on_level_size1}
s(2t) &\geq n \prod_{j=0}^{2t} q(j) - \frac{3}{2} \notag \\
&\geq 2^{2t} \left(\frac{1}{2}\right)^t \left(\frac{2}{5}\right)^{\frac{t}{3}}  \left(\frac{3}{5}\right)^{\frac{2 t}{3}+1} - \frac{3}{2} &&\text{as $|S'_R \bigcap S'_2| < \frac{t}{3}$} \notag \\
&\geq 4 &&\text{for large enough $t$}
\end{align}
Where the last inequality follows from the fact that $(\frac{2}{5}) \cdot (\frac{3}{5})^2 > (\frac{1}{2})^3$. \\
Now as $l \geq 2t$ and $|S'_R| < t$ we conclude that $|S'_B| \geq t$. Furthermore bound \ref{crude_bound_on_level_size1} implies that
for $0 \leq i \leq 2t$  the following holds for the bag sizes of RRT $F''$.
\begin{equation}\label{tedious_argument_but_we_need_it1}
s(i) \geq \frac{n}{2} \prod_{j=0}^{i} q(j)
\end{equation}
Assume that $|S'_B \bigcap S'_1| = \delta t$ for some $ 0 \leq  \delta < \frac{1}{3} $.
Let $Q'' \subseteq S'_B$ be a subset of $S'_B$ such that $|Q''|=t$ and $|Q'' \bigcap S'_1| = \delta t$.
Thus by Lemma \ref{tiny_lemma_bag_sizes10} (applied for RRT $F''$) we have
\begin{align}\label{onseventh_lemma_long_calc1}
M_t(G) &\geq 2^{-O(t \log t)}  \prod_{i \in Q''} s(i-1) \\
&\geq  2^{-O(t \log t)} \prod_{(1-\delta)t < i < t} s(i-1)  \prod_{(1+\delta)t < i < 2t} s(i-1) &&\text{as $|Q'' \bigcap S'_1| = \delta t$ }\notag \\
&\geq  2^{-O(t \log t)} \prod_{(1-\delta)t < i < t}  2^{2t-i} \prod_{0 < i < (1-\delta)t} 2^{2t}  \left(\frac{1}{2}\right)^t \left(\frac{2}{5}\right)^{\delta t}
\left(\frac{3}{5}\right)^{i} &&\text{by Inequality \ref{tedious_argument_but_we_need_it1}} \notag \\
&\geq 2^{-O(t \log t)}  \: 2^{( 1+\frac{1}{2}\delta^2 )t^2} \left(\frac{2}{5}\right)^{\delta(1-\delta)t^2} \left(\frac{3}{5}\right)^{\frac{1}{2}(1-\delta)^2 t^2}\notag \\
&\geq  2^{f(\delta)t^2-O(t \log t)} \notag
\end{align}
Where $f(\delta) =  1+\frac{1}{2}\delta^2 + \log(0.4)\delta(1-\delta) + \log(0.6)\frac{1}{2}(1-\delta)^2$.
Now one can verify that the quadratic polynomial $f(\delta)$ satisfies for all $\delta$ the following.
\begin{equation}\label{onseventh_lemma_long_calc2}
f(\delta) \geq \frac{4}{7} + 10^{-4}
\end{equation}
And we conclude from Inequalities \ref{onseventh_lemma_long_calc1} and \ref{onseventh_lemma_long_calc2} that for large enough $t$ we have
\[
M_t(G) \geq 2^{\frac{4}{7}t^2}
\]
This concludes the proof.
\EPF
\BCR
For large enough even $t$,  any $2$-coloring of the edges of the complete graph on $2^{t}$ vertices contains at least
$
2^{\frac{1}{7} t^2}
$
monochromatic complete subgraphs of size $t/2$.
\ECR

\section{Average size of a monochromatic complete subgraph}\label{average_section1}

Let $A(n)$ be the maximum number $m$ such that any $2$-coloring of the edges of a complete graph on $n$ vertices has an average size of monochromatic complete subgraphs at least $m$. The following theorem implies that $A(n) > 0.292\log n$, for sufficiently large $n$.

\BTHM
For large $t$ we have
\[
 A(2^t) \geq t\left( 1- \sqrt{\frac{1}{2}} -o(1) \right) > 0.292t
\]
\ETHM
\BPF
We will prove that for any fixed $\epsilon>0$ and large enough $t$ (depending on $\epsilon>0$),
the average size of a monochromatic complete subgraph in any $2$-coloring of a complete graph on $2^t$ vertices is at least
$t\left( 1- \sqrt{\frac{1}{2}} - 2\epsilon \right)$. \\
Fix $\epsilon>0$.
Let $G$ be a complete graph on $2^t$ vertices whose edges are colored with red and blue.
Let $F$ be the General Ramsey Tree associated with graph $G$.  \\
Set $r=\left\lceil t\left( 1- \sqrt{\frac{1}{2}} - \epsilon \right) \right\rceil$. Notice that for large enough $t$ we have $r \leq \frac{t}{3}$.
Let $Q_i$ denote the set of nodes on level $i$ of the GRT $F$.
Now by Lemma \ref{important_layer_bound_lem1} we have for all $0 \leq i \leq r$ and large enough $t$ that
\begin{equation}\label{Qbound1}
\frac{|Q_{i+1}|}{|Q_i|} \geq 2^{t-i} -2 \geq 2^{t-1-i}
\end{equation}
Hence
\begin{align}
\frac{|Q_{r+1}|}{|Q_0|} &= \prod_{i=0}^{r} \frac{|Q_{i+1}|}{|Q_i|} \label{product_ident1}\\ &\geq \prod_{i=0}^{r} 2^{t-1-i} \notag
\end{align}
Thus we may assume that
\begin{equation}\label{average_exact_bound1}
\frac{|Q_{r+1}|}{|Q_0|} =  2^{ \frac{\delta}{4} t^2}  \prod_{i=0}^{r} 2^{t-1-i}
\end{equation}
for some $\delta \geq 0$.
We conclude from Equations \ref{product_ident1} and  \ref{average_exact_bound1} that
\[
\prod_{i=0}^{r} \frac{|Q_{i+1}|}{|Q_i|} =  \prod_{i=0}^{r} 2^{t-1-i + \frac{\delta}{4}t^2 \frac{1}{r+1} }
\]
And hence for some $0 \leq j \leq r$
we have
\begin{align}\label{one_more_bound1}
\frac{|Q_{j+1}|}{|Q_j|} &\geq 2^{t-1-j + \frac{\delta}{4}t^2 \frac{1}{r+1} }  \notag  \\
&\geq 2^{t(1+\frac{3\delta}{4})-1-j} &&\text{as $r+1 \leq \frac{t}{3}$ for large enough $t$.}
\end{align}
Now by Lemma \ref{important_layer_bound_lem2} and Inequality \ref{one_more_bound1}
we have that
for all $j \leq i \leq t-3$ the following holds.
\begin{align}\label{Qbound2}
\frac{|Q_{i+1}|}{|Q_i|} &\geq 2^{t(1+\frac{3\delta}{4})-1-i} - 2 \\
&\geq 2^{t(1+\frac{3\delta}{4})-2-i} \notag
\end{align}
We conclude from bounds \ref{Qbound1} and \ref{Qbound2} that
\begin{align}
  \frac{|Q_{t-2}|}{|Q_0|} &= \prod_{i=0}^{t-3} \frac{|Q_{i+1}|}{|Q_i|} \notag \\
  &=  \prod_{i=0}^{j-1} \frac{|Q_{i+1}|}{|Q_i|}  \cdot  \prod_{i=j}^{t-3} \frac{|Q_{i+1}|}{|Q_i|} \notag \\
  &\geq \prod_{i=0}^{j-1} 2^{t-1-i} \cdot \prod_{i=j}^{t-3}  2^{t(1+\frac{3\delta}{4})-2-i} \notag \\
  &\geq \prod_{i=0}^{t-3} 2^{t-i} \cdot \prod_{i=r}^{t-3} 2^{t(\frac{3\delta}{4})} \cdot 2^{-O(t)} &&\text{as $j \leq r$} \notag \\
  &\geq 2^{\frac{1}{2}t^2} \cdot 2^{(t-r)t(\frac{3\delta}{4})} \cdot 2^{-O(t)} \notag \\
  &\geq 2^{\frac{1}{2}(1 + \delta)t^2 - O(t)} &&\text{as $r \leq \frac{t}{3}$ for large $t$} \notag
\end{align}
Thus
\begin{equation}\label{lower_clique_bound_average1}
  |Q_{t-2}| \geq  2^{\frac{1}{2}(1 + \delta)t^2 - O(t)}
\end{equation}
And we conclude by Lemma \ref{anne_idea1_average_section} and Inequality \ref{lower_clique_bound_average1} that $G$ contains at least $ 2^{\frac{1}{4}(1 + \delta)t^2 - O(t \log t)}  $
monochromatic complete subgraphs.
Now we shall give an upper-bound on the number of monochromatic complete subgraphs of size at most $r+2$ in graph $G$.
By Equation \ref{average_exact_bound1} we have
\begin{align}\label{yet_another_horrible_calculation1}
\frac{|Q_{r+1}|}{|Q_0|} &=  2^{ \frac{\delta}{4} t^2}  \prod_{i=0}^{r} 2^{t-1-i} \notag \\
&\leq  2^{\frac{\delta}{4} t^2 + t r - \frac{r^2}{2} +O(t)} \notag \\
&\leq  2^{\frac{t^2}{4}(1 + \delta - 2\epsilon^2 - 2\sqrt{2}\epsilon ) +O(t)}
\end{align}
Hence by Inequality \ref{yet_another_horrible_calculation1} and the fact that $|Q_0|=2^t$ we have
\begin{equation}\label{upper_bound_on_average_size11}
  |Q_{r+1}| \leq 2^{\frac{t^2}{4}(1 + \delta - \epsilon ) +O(t)}
\end{equation}
Furthermore by Equation \ref{Qbound1} we have for all $i \leq r$ (and large enough $t$) that $|Q_i| \leq |Q_{r+1}|$.
Thus by Lemma \ref{anne_idea1_average_section2} and Inequality \ref{upper_bound_on_average_size11} graph $G$ contains at most
\[
2^{\frac{t^2}{4}(1 + \delta - \epsilon ) +O(t)}
\]
monochromatic complete subgraphs of size at most $r+2 = \left\lceil t\left( 1- \sqrt{\frac{1}{2}} - \epsilon \right) \right\rceil + 2$.
Hence we have proven the following
\begin{itemize}
  \item Graph $G$ contains at least $ 2^{\frac{t^2}{4}(1 + \delta) - O(t \log t)}$  monochromatic complete subgraphs.
  \item Graph $G$ contains at most $2^{\frac{t^2}{4}(1 + \delta - \epsilon) +O(t)}$ monochromatic complete subgraphs of size at most
  $t\left( 1- \sqrt{\frac{1}{2}} - \epsilon \right)$.
\end{itemize}
We conclude that the average size of a monochromatic complete subgraph in $G$ is at least $t\left( 1- \sqrt{\frac{1}{2}} - 2\epsilon \right)$
for large enough $t$ (depending on $\epsilon$) and thus we are done.
\EPF

Random graphs provide an upper bound on $A(n)$.

\BTHM \label{average_at_most_log1}
There are graphs on $n$ vertices for which the average size of monochromatic complete subgraphs satisfies $A(n) \leq (1+o(1))\log n$.
\ETHM

\BPF
Let $G$ be a complete graph on $n$ vertices in which each edge is colored red with probability $\frac{1}{2}$ and blue otherwise.
Let $C$ be the number of monochromatic complete subgraphs in $G$.
In the Corollary of Theorem $3.1$ in \cite{olejar1997order} it is shown that with probability $1-o(1)$, namely asymptotically always surely (a.a.s.), we have for large enough $n$  the following
\begin{equation}\label{olejar_thm1}
C \geq n^{\frac{1}{2}\log n - \log \log n} = 2^{\frac{1}{2}(\log n)^2 - \log n \log \log n}
\end{equation}
Let $C_t$ be the number of monochromatic complete subgraphs of size $t$  in $G$. \\
Let $C_{\geq t}$ be the number of monochromatic complete subgraphs of size at least $t$  in $G$.
Notice that
\begin{align}\label{expected_val_avg1}
E(C_t) &= \binom{n}{t} 2^{-\binom{t}{2} +1} \\ & \leq n^t 2^{-\frac{1}{2}t^2} &&\text{for large enough $t$} \notag
\\ &= 2^{t(\log n -\frac{1}{2}t)} \notag
\end{align}
Hence given a fixed $\epsilon>0$ we have
\begin{align}\label{expected_val_avg2}
  E\left(C_{\geq (1+\epsilon)\log n }\right) &\leq n 2^{(1+\epsilon)\log n(\log n -\frac{1}{2}(1+\epsilon)\log n)} \\
  &\leq 2^{\frac{1}{2}(1-\epsilon)(\log n)^2} &&\text{for large enough $n$} \notag
\end{align}
Thus by Markov's inequality we have a.a.s. that
\begin{equation}\label{expected_val_avg3}
 C_{\geq (1+\epsilon)\log n } \leq n 2^{\frac{1}{2}(1-\epsilon)(\log n)^2}
\end{equation}
Hence a.a.s. Graph $G$ satisfies both inequalities \ref{expected_val_avg3} and \ref{olejar_thm1}. \\
Let $C' =  C_{\geq (1+\epsilon)\log n }$.
We conclude that
\begin{align}
A(n) &\leq \frac{    (C - C' ) (1+\epsilon)\log n  + C' n }{C} \notag \\
&\leq (1+\epsilon)\log n + n\frac{C'}{C} \notag \\
&\leq (1+2\epsilon)\log n &&\text{for large enough $n$} \notag
\end{align}
Thus we are done.
\EPF

\section{Relations between the different parameters}\label{parm_relationship1}

\subsection{Relations between the average size of a monochromatic complete subgraph and the number of monochromatic complete subgraphs}

In this section we will show certain relations between the average size of a monochromatic complete subgraph and the number of monochromatic complete subgraphs in a $2$-coloring of a complete graph.
Let $A(n)$ be the maximum number $m$ such that any $2$-coloring of the edges of a complete graph on $n$ vertices
has an average size of monochromatic complete subgraphs at least $m$.
Let $C(n)$ be the maximum number $m$ such that any $2$-coloring of the edges of a complete graph on $n$ vertices
contains at least $m$ monochromatic complete subgraphs.
\BTHM\label{average_to_number_of_monocliques1}
For all $n \geq 2$, we have
\[
C(n) \geq e^{\sum_{i=2}^{n} \frac{A(i)}{i}}
\]
\ETHM
\BPF
Let $G$ be an arbitrary $2$-coloring of the complete graph on $n$ vertices.
Suppose that the vertices of $G$ are $v_1,v_2,\ldots,v_n$.
Denote by $A(G)$ the average size of a monochromatic complete subgraph in $G$.
Denote by $C(G)$ the number of monochromatic complete subgraph in $G$.
Denote by $C_i(G)$ the number of monochromatic complete subgraph in $G$ that contain vertex $v_i$.
Denote by $G - v_j$ the graph $G$ from which we remove the vertex $v_j$.
Notice that
\begin{equation}\label{avg_total_eq1}
\frac{\sum_{i=1}^{n} C_i(G)}{C(G)} = A(G)
\end{equation}
Furthermore as there are $n$ vertices in $G$ we have for some index $j$ that
\begin{equation}\label{avg_total_eq2}
C_j(G) \geq \frac{\sum_{i=1}^{n} C_i(G)}{n}
\end{equation}
Thus we have the following for some index $j$
\begin{align}\label{avg_total_eq3}
C(G) &= C_j(G) + C(G-v_j) \\
&\geq  \frac{\sum_{i=1}^{n} C_i(G)}{n} + C(G-v_j) &&\text{by Inequality \ref{avg_total_eq2}} \notag \\
&\geq \frac{1}{n} C(G)A(G) + C(G-v_j)  &&\text{by Equation \ref{avg_total_eq1}} \notag \\
&\geq \frac{1}{n} C(n)A(n) + C(n-1) \notag
\end{align}
As Inequality \ref{avg_total_eq3} holds for any graph $G$ on $n$ vertices we conclude that
\begin{equation}\label{avg_total_eq4}
C(n) \geq \frac{1}{n} C(n)A(n) + C(n-1)
\end{equation}
Hence
\begin{equation}\label{avg_total_eq5}
C(n)\left(1 -  \frac{A(n)}{n}  \right) \geq  C(n-1)
\end{equation}
Now as  $1-\frac{A(n)}{n} \leq e^{-\frac{A(n)}{n}}$ we conclude that
\begin{align}\label{avg_total_eq6}
  C(n) &\geq e^{\frac{A(n)}{n}} C(n-1)
\end{align}
Iterating the argument we get
\[
C(n) \geq e^{\sum_{i=2}^{n} \frac{A(i)}{i}}
\]
and thus we are done.
\EPF
\BCR\label{conditional_cliques_on_avg_cor1}
If there is some constant $c$ such that for all large enough $n$, $A(n) \geq c\log n$, then for all large enough $n$, $C(n) = \Omega\left(  n^{\frac{c}{2}\log n}\right)$.
\ECR
\BPF
If for all $n \geq n_0$ we have $A(n) \geq c\log n$, then by Theorem \ref{average_to_number_of_monocliques1} we have the following for all $n \geq n_0$.
\begin{align}\label{average_to_number_of_monocliques2}
C(n)  &\geq e^{\sum_{i=2}^{n} \frac{A(i)}{i}} \\
&\geq e^{\sum_{i=n_0}^{n} \frac{c\log i}{i}} \notag \\
&\geq e^{\int_{n_0}^{n} \frac{c \log x}{x} dx} \notag \\
&= \Omega\left(  n^{\frac{c}{2}\log n}\right) &&\text{as $\int  \frac{c \log x}{x} dx = \frac{c}{2} \log(x) \ln(x)$}\notag
\end{align}
And thus we are done.
\EPF
We note that for $c=1$ Corollary \ref{conditional_cliques_on_avg_cor1} is tight for the random graph $G(n,\frac{1}{2})$, as $A(G)$ is roughly $\log n$, while $C(G)$ is roughly $n^{\frac{1}{2}\log n}$.

\subsection{Relations between the average size of a monochromatic complete subgraph and the maximum size of a monochromatic complete subgraph }

To simplify the notation in all this section when we say a graph $G$ on $n$ vertices we mean that $G$ is a $2$-coloring of the complete graph on $n$ vertices.
Denote by $A(G)$ the average size of a monochromatic complete subgraph in $G$. Denote by $M(G)$ the maximum size of a monochromatic complete subgraph in $G$.

\BTHM
For large enough $n$, any graph $G$ on $n$ vertices satisfies $\frac{M(G)}{A(G)} \leq (1+o(1))\log n$.
\ETHM

\BPF
As $A(G) = \Omega(\log n)$ we may assume that $M(G) = \Omega(\log^2 n)$.
Suppose that $M(G)=k$. Hence graph $G$ contains at least $2^k$ monochromatic complete subgraphs.
Now notice that for any fixed constant $\epsilon>0$, the number of monochromatic complete subgraphs in $G$ of size at most $(1-\epsilon)\frac{k}{\log n}$
is at most
\begin{equation}\label{Anne_max_avg7}
  \sum_{i=0}^{\frac{(1-\epsilon) k}{\log n}} \binom{n}{i} \leq k \binom{n}{\frac{(1-\epsilon) k}{\log n}} \leq k 2^{(1-\epsilon) k} \notag
\end{equation}
We conclude that
\begin{align}
A(G) &\geq \frac{(2^k - k 2^{(1-\epsilon) k}) \frac{(1-\epsilon) k}{\log n}}{2^k} \notag \\
 &\geq \frac{(1-2\epsilon) k}{\log n} &&\text{for large enough $n$} \notag
\end{align}
Where the last inequality follows from the fact that for large enough $n$ we have
\[
\frac{2^k - k 2^{(1-\epsilon) k}}{2^k} \geq 1-\epsilon
\]
And thus we are done.
\EPF

\BTHM
For infinitely many values of $n$, there is graph $G$ on $n$ vertices which satisfies $\frac{M(G)}{A(G)} \geq \left(\frac{1}{2}-o(1)\right)\log n$ and $M(G) = O(\log^2 n)$.
\ETHM
\BPF
We recall our definitions.
Denote by $C(G)$  the number of monochromatic complete subgraphs in graph $G$.
Denote by $A(G)$  the average size of a monochromatic complete subgraph in graph $G$. \\
Denote by $C_b(G)$  the number of blue monochromatic complete subgraphs in graph $G$.
Denote by $A_b(G)$  the average size of a blue monochromatic complete subgraph in graph $G$. \\
Denote by $C_r(G)$  the number of red monochromatic complete subgraphs in graph $G$.
Denote by $A_r(G)$  the average size of a red monochromatic complete subgraph in graph $G$. \\
Denote by $M(G)$ the maximum size of a monochromatic complete subgraph in graph $G$. \\
Let $G_1$ be a complete graph on $t$ vertices in which each edge is colored red with probability $\frac{1}{2}$ and blue otherwise. \\
Using the same techniques as in Theorem \ref{average_at_most_log1} we can assume that graph $G_1$ satisfies a.a.s. the following properties
for large $t$.
\begin{enumerate}
  \item $C_b(G_1) \geq 2^{\frac{1}{2}\log^2 t - \log t \log \log t}$
  \item $A_b(G_1) \leq (1+o(1))\log t$
  \item  $C_r(G_1) \geq 2^{\frac{1}{2}\log^2 t - \log t \log \log t}$
  \item $A_r(G_1) \leq (1+o(1))\log t$
\end{enumerate}
Fix an arbitrary small $\epsilon > 0$.
Let $G_2$ be a red monochromatic complete subgraph on $\left(\frac{1}{2}-\epsilon\right)\log^2 t$ vertices.
Let graph $G$ be the union of graphs $G_1$ and $G_2$ where each edge between the two graphs is colored blue.
Without loss of generality assume that the color in $G$ of the empty set and the sets of size $1$ is blue.
By the linearity of expectation of the average operator
and the fact that $G_2$ is a red monochromatic complete subgraph we have
\begin{equation}\label{average_magic1_lem1}
  A_b(G) \leq A_b(G_1) + 1 \leq (1+o(1))\log t
\end{equation}
Furthermore as $C_r(G) = C_r(G_1) + C_r(G_2)$, we have $A_r(G) \leq A_r(G_1) + \frac{n C_r(G_2)}{C_r(G_1)}$.
And hence as  $C_r(G_1) \geq 2^{\frac{1}{2}\log^2 t - \log t \log \log t}$ while $C_r(G_2) \leq 2^{\left(\frac{1}{2}-\epsilon\right)\log^2 t}$ and
as we can take arbitrarily small $\epsilon$ and large enough $t$ (depending on $\epsilon$), we get that
\begin{equation}\label{average_magic1_lem2}
  A_r(G) \leq A_r(G_1) +o(1) = (1+o(1))\log t
\end{equation}
Now as $A(G) \leq \max\{A_r(G),A_b(G)\}$ we conclude from Inequalities \ref{average_magic1_lem1} and \ref{average_magic1_lem2} that
\begin{equation}\label{final_estimate_on_avgofg}
A(G) \leq (1+o(1))\log t
\end{equation}
We also notice that
\begin{equation}\label{final_estimate_on_max}
M(G) \geq M(G_2) = \left(\frac{1}{2}-o(1)\right)\log^2 t
\end{equation}
We conclude from inequalities \ref{final_estimate_on_avgofg} and \ref{final_estimate_on_max} that
$\frac{M(G)}{A(G)} \geq \left(\frac{1}{2}-o(1) \right)\log t$
where $G$ is a graph on $n=t+ \left(\frac{1}{2}-o(1)\right)\log^2 t$ and the theorem follows as $\log(t) \geq \log(n)-1$.
\EPF
\BTHM
For infinitely many values of $n$, there is graph $G$ on $n$ vertices which satisfies $M(G) -A(G) \leq 1 + o(1)$ and $M(G) = O(\log^2 n)$.
\ETHM
\BPF
We recall our definitions.
Denote by $C(G)$  the number of monochromatic complete subgraphs in graph $G$.
Denote by $A(G)$  the average size of a monochromatic complete subgraph in graph $G$. \\
Denote by $C_b(G)$  the number of blue monochromatic complete subgraphs in graph $G$.
Denote by $A_b(G)$  the average size of a blue monochromatic complete subgraph in graph $G$. \\
Denote by $C_r(G)$  the number of red monochromatic complete subgraphs in graph $G$.
Denote by $A_r(G)$  the average size of a red monochromatic complete subgraph in graph $G$. \\
Denote by $M(G)$ the maximum size of a monochromatic complete subgraph in graph $G$. \\
Let $G_1$ be a complete graph on $2^{\sqrt{t}}$ vertices in which each edge is colored red with probability $\frac{1}{2}$ and blue otherwise. \\
From standard estimates on random graphs we can assume that graph $G_1$ satisfies a.a.s. the following properties
for large $t$.
\begin{enumerate}
  \item $M(G_1) < 3\sqrt{t}$
  \item $C(G_1) \leq 2^{O(t)}$
\end{enumerate}
Now Let $G_2$ be a complete graph on $t(t-3\sqrt{t})$ vertices consisting of $t$ red monochromatic complete subgraphs $R_1,R_2,\ldots,R_t$, where for all $1\leq i \leq t$ we have
$|R_i|= t-3\sqrt{t}$, furthermore each edge between two vertices in different $R_i$'s is colored blue. \\
Graph $G$ is the union of graphs $G_1$ and $G_2$ where each edge between the two graphs is colored red.
Without loss of generality assume that the color in $G$ of the empty set and the sets of size $1$ is blue.
By the linearity of expectation of the average operator we have
\begin{align}\label{lin_of_expectation_argument2}
  A_b(G_2) &= \sum_{i=1}^{t} A_b(R_i) \\ &= t \cdot \frac{t-3\sqrt{t}}{t-3\sqrt{t}+1}
  \geq t - 1 - o(1) \notag
\end{align}
Also notice that
\begin{equation}\label{cb_size_estimate_1}
  C_b(G_2) \geq (t-3\sqrt{t})^t \geq 2^{t\log t-O(t)}
\end{equation}
Furthermore we have
\begin{align}\label{main_size_estimate_avg_lem2}
  C(G) - C_b(G_2) &\leq (C(G_1)+1) \cdot (C_r(G_2)+1) \\
  &\leq 2^{O(t)}  \notag
\end{align}
Where the last inequality follows from the fact that $C_r(G_2) \leq  2^{O(t)}$ and $C(G_1) \leq  2^{O(t)}$.
Hence we have
\begin{align}\label{tying_up_avg_lower_bound1}
  A(G) &\geq A_b(G_2) \cdot \frac{C_b(G_2)}{C(G)} \\
  &= A_b(G_2) - A_b(G_2) \cdot \frac{C(G) - C_b(G_2)}{C(G)} \notag \\
  &\geq A_b(G_2) - A_b(G_2) \cdot \frac{C(G) - C_b(G_2)}{C_b(G_2)} \notag \\
  &\geq A_b(G) - o(1) \notag
\end{align}
Where the last inequality follows from Inequalities \ref{cb_size_estimate_1} and \ref{main_size_estimate_avg_lem2}.
Hence from Inequalities \ref{tying_up_avg_lower_bound1} and \ref{lin_of_expectation_argument2} we have
\[
A(G) \geq A_b(G) - o(1) \geq t-1-o(1)
\]
Furthermore we have $M(G) = t$, as $M(G_2)=t$ and $M(G_1) < 3\sqrt{t}$.
Now as $G$ is a graph on $n=2^{\sqrt{t}} + t(t-3\sqrt{t})$ vertices, the theorem follows.
\EPF

\BTHM
For every integer $c\geq 0$ the following holds. For infinitely many values of $n$, there is graph $G$ on $n$ vertices which
satisfies $\frac{M(G)}{A(G)} \leq 1 + \frac{1}{c+1} + o(1)$ and $M(G) \leq (c+2+o(1))\log n$.
\ETHM
\BPF
We recall our definitions.
Denote by $C(G)$  the number of monochromatic complete subgraphs in graph $G$.
Denote by $A(G)$  the average size of a monochromatic complete subgraph in graph $G$. \\
Denote by $C_b(G)$  the number of blue monochromatic complete subgraphs in graph $G$.
Denote by $A_b(G)$  the average size of a blue monochromatic complete subgraph in graph $G$. \\
Denote by $C_r(G)$  the number of red monochromatic complete subgraphs in graph $G$.
Denote by $A_r(G)$  the average size of a red monochromatic complete subgraph in graph $G$. \\
Denote by $M(G)$ the maximum size of a monochromatic complete subgraph in graph $G$. \\
Let $G_1$ be a complete graph on $2^{t}$ vertices in which each edge is colored red with probability $\frac{1}{2}$ and blue otherwise. \\
We assume without loss of generality the $C_b(G_1) \geq C_r(G_1)$.
From standard estimates on random graphs and techniques similar to those in Theorem \ref{average_at_most_log1}
we can assume that graph $G_1$ satisfies a.a.s. the following properties
for large $t$.
\begin{enumerate}
  \item $M(G_1) \leq  (2+o(1))t$
  \item $A_b(G_1) \geq  (1-o(1))t$
\end{enumerate}
Now Let $G_2$ be a complete graph on $(c t)^2$ vertices consisting of $c t$ red monochromatic complete subgraphs $R_1,R_2,\ldots,R_{c t}$, where for all $1\leq i \leq c t$ we have
$|R_i|= c t$, furthermore each edge between two vertices in different $R_i$'s is colored blue. \\
Graph $G$ is the union of graphs $G_1$ and $G_2$ where each edge between the two graphs is colored blue.
Without loss of generality assume that the color in $G$ of the empty set and the sets of size $1$ is blue.
By the linearity of expectation of the average operator we have
\begin{align}
  A_b(G_2) &= \sum_{i=1}^{t} A_b(R_i) \\ &= ct \cdot \frac{ct}{ct+1}
  \geq ct - 1 \notag
\end{align}
Applying linearity of expectation once again we get
\begin{equation}\label{blue_aerage1_final_lem}
  A_b(G) = A_b(G_1) + A_b(G_2) \geq (c+1-o(1))t
\end{equation}
Now
\[
A(G) \geq A_b(G) \cdot \frac{C_b(G)}{C(G)} \geq A_b(G) (1-o(1))
\]
Furthermore notice that $M(G) \leq (c+2+o(1))t$ and thus
\[
\frac{M(G)}{A(G)} \leq 1 + \frac{1}{c+1} + o(1)
\]
and we are done.
\EPF

\BTHM
\label{thm:mode}
Any graph $G$ on $n\geq 6$ vertices satisfies $M(G) - A(G) \geq \frac{1}{2}$.
\ETHM

\BPF
Let $M$ be the size of a maximum monochromatic complete subgraph in $G$.
As $n\geq 6$ we have $M \geq 3$.
Let $C_i$ be the number of monochromatic complete subgraphs of size $i$ in $G$.
We will prove that $C_M \leq C_{M-1}$. \\
Notice that each monochromatic complete subgraph of size $M$ in $G$ contains $M$ monochromatic complete subgraphs of size $M-1$.
Hence there exist at least one monochromatic complete subgraph of $C$ of size $M-1$ that is a subgraph of $\frac{M C_M}{C_{M-1}}$ different monochromatic complete subgraphs of size $M$. Without loss of generality, let $C$ be a red monochromatic complete subgraph. Taking all the vertices from the size $M$ monochromatic complete subgraphs that are supersets of $C$, excluding the vertices of $C$, we have a subgraph of at least $\frac{M C_M}{C_{M-1}}$ vertices.
This subgraph cannot contain a red edge, or else the edge's two endpoint vertices
and $C$ would form a size $M + 1$ red monochromatic complete subgraph in $G$.
This subgraph also cannot contain a size $M + 1$ blue monochromatic complete subgraph. Therefore, we have the following inequality:
\[
\frac{M C_M}{C_{M-1}} \leq M
\]
And we conclude that $C_M \leq C_{M-1}$. Hence
\[
A(G) \leq \frac{M C_M + (M-1)C_{M-1}}{C_M + C_{M-1}} \leq M - \frac{1}{2}
\]
And we are done.
\EPF

Using the notation of Theorem~\ref{thm:mode}, the proof of Theorem~\ref{thm:mode} shows that $C_M \le C_{M-1}$. Viewing a single vertex as a complete subgraph in each of the two colors, the 5-cycle is an example where $M=2$ and $C_2 = C_1 = 10$. The Paley graph on 17 vertices (the 5-cycle is a Paley graph on 5 vertices) has 136 edges, 136 triangles, and no clique of size~4. Being self complimentary, for it we have that $M=3$ and $C_3 = C_2 = 272$.

\section*{Acknowledgements}
Work supported in part by the Israel Science Foundation (grant No. 1388/16).
An earlier partial version of this work  was posted on the arXiv (https://arxiv.org/abs/1703.09682) in March 2017. We thank David Conlon for bringing to our attention the work of Sz{\'{e}}kely.

\bibliographystyle{alpha}

\begin{appendix}
\section{ A few facts which are needed for Theorem \ref{cool_ramsey_theorem1} }\label{first_appendix}
\BL\label{tiny_cool_bound1}
Given an integer $t\geq 1$ we have
\[
\sum_{ \substack{ S \subset \mathbb Z_{\ge 0} \\ |S|=t}} 2^{-  \sum_{i \in S} i} = \frac{2^t}{\prod_{i=1}^{t} {(2^i-1)}  }
\]
\EL
\BPF
We shall prove the claim by induction on $t$.
For $t=1$ we have
\[
\sum_{ \substack{ S \subset \mathbb Z_{\ge 0} \\ |S|=1}} 2^{-  \sum_{i \in S} i} = \sum_{i\geq 0} 2^{-i} = 2
\]
Now we assume that the claim holds for $t-1$ and will prove for $t$.
Recall that $Z_{\ge j}$ denotes the set of all integers greater or equal to $j$.
\begin{align}\label{tiny_cool_bound1_calculation1}
\sum_{ \substack{ S \subset \mathbb Z_{\ge 0} \\ |S|=t}} 2^{-  \sum_{i \in S} i} &=
\sum_{j \geq 0} \sum_{ \substack{ S \subset \mathbb Z_{\ge j+1} \\ |S|=t-1}} 2^{-j -  \sum_{i \in S} i} \notag \\
&=  \sum_{j \geq 0} \sum_{ \substack{ S \subset \mathbb Z_{\ge 0} \\ |S|=t-1}} 2^{-j -(j+1)(t-1) -  \sum_{i \in S} i} \notag \\
&= \sum_{j \geq 0} 2^{-j -(j+1)(t-1)} \sum_{ \substack{ S \subset \mathbb Z_{\ge 0} \\ |S|=t-1}} 2^{ -  \sum_{i \in S} i} \notag \\
&= \frac{2}{2^t - 1}\sum_{ \substack{ S \subset \mathbb Z_{\ge 0} \\ |S|=t-1}} 2^{-  \sum_{i \in S} i} \notag \\
&= \frac{2^t}{\prod_{i=1}^{t} {(2^i-1)}} \notag
\end{align}
where the last equality follows from the induction hypothesis. Thus we are done.
\EPF
\BL\label{tiny_cool_bound2}
\[
\prod_{i=1}^{\infty} \Big(1-\frac{1}{2^i}\Big) \geq \frac{1}{4}
\]
\EL
\BPF
The pentagonal number theorem (proven by Euler) states that for all real $|x|<1$ we have
\begin{align}
\prod_{i=1}^{\infty} (1-x^i) &= 1 + \sum_{i=1}^{\infty} (-1)^i (x^{i(3i+1)/2} + x^{i(3i-1)/2}  ) \notag \\
&= 1 - x- x^2 +x^5 + x^7 - x^{12} - x^{15} + \ldots \notag
\end{align}
Hence we have for all real $0\leq x<1$
\begin{equation}\label{pentagonal2}
\prod_{i=1}^{\infty} (1-x^i) \geq 1 - x - x^2
\end{equation}
We set $x=\frac{1}{2}$ in Equation \ref{pentagonal2} and we are done.
\EPF

\end{appendix}
\end{document}